\documentclass[preprint,12pt,authoryear]{elsarticle}

\usepackage{color}
\usepackage{ulem}
\usepackage{amsmath,amssymb}
\usepackage{graphicx,psfrag,epsf}
\usepackage{enumerate}
\usepackage{natbib}
\usepackage{url}
\usepackage{makecell}
\usepackage{booktabs,tabularx}

\DeclareRobustCommand{\rpsi}{\text{\reflectbox{$\psi$}}}

\journal{Alexandria Engineering Journal}

\begin{document}

\begin{frontmatter}



\title{Uncovering a generalised gamma distribution: from shape to interpretation}

\author[1]{Matthias Wagener}
\author[1]{Andriette Bekker}
\author[2]{Mohammad Arashi}
\author[3]{Antonio Punzo}

\affiliation[1]{organization={Department of Statistics, University of Pretoria}, addressline={Lynnwood Street}, city={Pretoria}, postcode={0002}, state={Gauteng}, country={South Africa}}
\affiliation[2]{organization={Department of Statistics, Ferdowsi University of Mashhad}, addressline={Azadi Square}, city={Mashhad}, postcode={9177948974}, state={Razavi Khorasan}, country={Iran}}
\affiliation[3]{organization={Department of Statistics, University of Catania}, addressline={Piazza Università}, city={Catania}, postcode={95124 }, state={Catania}, country={italy}}




\begin{abstract}
In this paper, we introduce the flexible interpretable gamma (FIG) distribution which has been derived by Weibullisation of the body-tail generalised normal distribution. The parameters of the FIG have been verified graphically and mathematically as having interpretable roles in controlling the left-tail, body, and right-tail shape. The generalised gamma (GG) distribution has become a staple model for positive data in statistics due to its interpretable parameters and tractable equations. Although there are many generalised forms of the GG which can provide better fit to data, none of them extend the GG so that the parameters are interpretable. Additionally, we present some mathematical characteristics and prove the identifiability of the FIG parameters. Finally, we apply the FIG model to hand grip strength and insurance loss data to assess its flexibility relative to existing models.
\end{abstract}



\begin{keyword}
body shape \sep heavy-tails \sep positive data \sep insurance losses \sep interpretability \sep maximum likelihood estimation
\end{keyword}

\end{frontmatter}

\section{Introduction}
\label{sec:intro}
\begin{figure}[h!]
  \includegraphics[width=13cm,keepaspectratio]{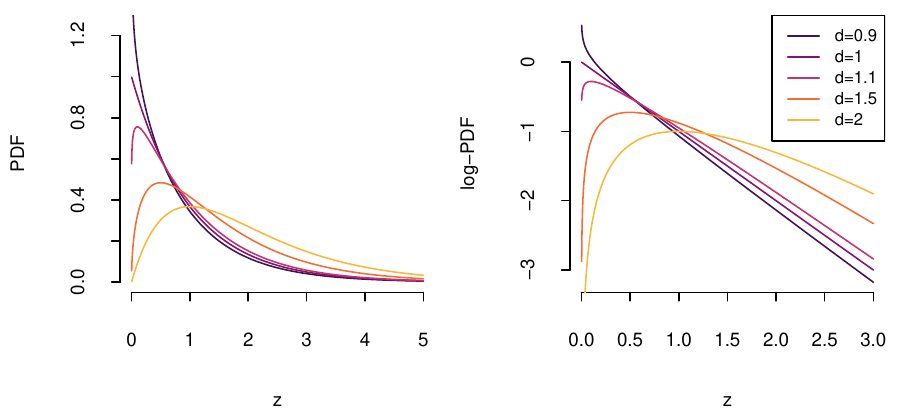}
 \caption{The GG PDF for different values of left-tail shape parameter $d$ and fixed right-tail shape $p=1$.}
 \label{fig:ggden_nu}
\end{figure}

\begin{figure}[h!]
  \includegraphics[width=13cm,keepaspectratio]{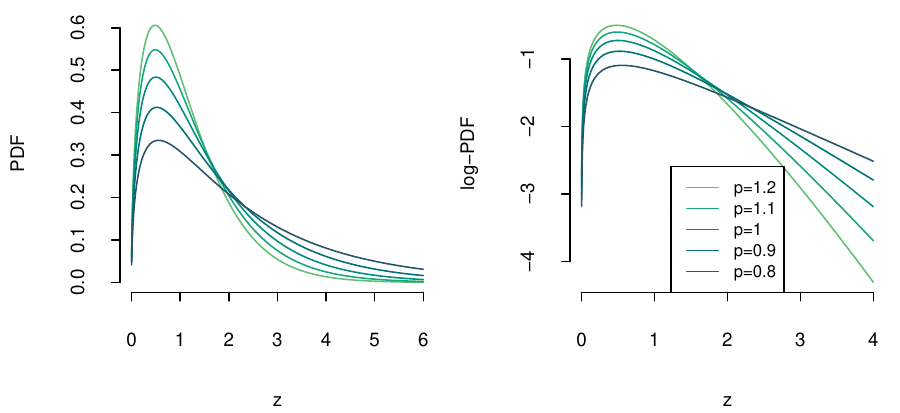}
 \caption{The GG PDF for different values of the right-tail shape parameter $p$ and fixed left-tail shape $d=1.5$.}
 \label{fig:ggden_beta}
\end{figure}

%
%
%
The best known form of the generalised gamma (GG) distribution was defined by Stacy in 1962 \cite{ggStacy}. Before this, a precursor model had been analysed by Amoroso in 1925 for income distribution modelling purposes \cite{Amoroso}. The GG contains numerous sub-models, including the exponential, gamma, Weibull, and log-normal, as limiting cases. The probability PDF of the standard (unscaled) GG is given below:
\begin{equation}
\label{eq:ggden}
f(z;p,d) = \frac{p}{\Gamma(d/ p)} z^{d-1}\text{e}^{-z^p},
\end{equation}
where $z,d,p>0$, $\Gamma(\cdot)$ is the gamma function, and is denoted as $Z\sim GG(p,d)$. Note that (\ref{eq:ggden}) is equivalent to having $a=1$ in the PDF of \cite{ggStacy}. This is done as a simplification, with the knowledge that a simple scaling can be applied after generalisation. The role of the GG distribution shape parameters become apparent when considering the derivative of the log of the kernel in (\ref{eq:ggden}),
\begin{equation}
\label{eq:ggkder}
d(z;p,d)=\frac{\partial}{\partial z}\ln{\left(z^{d-1}e^{-z^p}\right)}
=\frac{d-1}{z}-p z^{p-1}.
\end{equation}
    The left tail behaviour for GG is determined by $d$. Considering the case where, $d \neq 1$ and $\lim_{z\to 0^+}d(z;p,d)$ in (\ref{eq:ggkder}), the first term dominates as the second approaches zero making $d$ the primary shape determinant. If $d=1$, it has no affect on the left tail shape. Figure \ref{fig:ggden_nu} illustrates these GG PDF properties in relation to $d$. The right tail behaviour for the GG is influenced by $p$. Considering the case where, $\lim_{z\to\infty}d(z;p,d)$ in (\ref{eq:ggkder}), the first term approaches zero while the second dominates, thus making $p$ the primary shape determinant. Figure \ref{fig:ggden_beta} illustrates these GG PDF properties in relation to $p$. The broad range of distribution shapes for the GG can model has led to widespread application, such as in survival analysis, \cite{GGleegros}, time series \cite{GGARMA}, phonemic segmentation \cite{GGTextapp}, wireless fading models \cite{GGGwirelessapp}, drought \cite{GGNadarajahapp}, statistical size \cite{kleiber2003statistical}, demographic research \cite{CoaleNuptiDemo}, and economics \cite{Amoroso}.
\begin{table}[h!]
\footnotesize
    \centering
\begin{tabularx}{\textwidth} {
  | >{\centering\arraybackslash}X 
  | >{\centering\arraybackslash}X 
  | >{\centering\arraybackslash}X
  | >{\centering\arraybackslash}X |}
     \hline
     Distribution  & Number of parameters & Author/s & Year\\
     \hline
    GG &$3$ &  \cite{Amoroso}, \cite{ggStacy} & 1925,1962\\
     quotients of GG  &$3\text{ or }9ik \text{ for } i,k=1,2,3,\dots$ & \cite{QuoGG}, \cite{2020IGG}, \cite{RGGBilankulu} & 1967,2020,2021\\
    log-GG &$3$ &  \cite{GGCosulJain}, \cite{prenticelogGG}, \cite{bellreparGG} & 1971,1974,1988 \\
    unified GG &$4$ &  \cite{IGGagarkall}& 1996\\
    composite Weibullised GG & $5$ & \cite{pauw2010densities} & 2010\\
    exponentiated GG &$4$  & \cite{EGG} & 2011\\
    Kumaraswamy GG & $5$ & \cite{KGG} & 2011\\
    beta GG &$5$  & \cite{BGG}& 2012\\
    Kummer beta GG &$6$  & \cite{KBGG} & 2014\\
    transmuted GG & $4$ & \cite{transmGG}& 2015\\%
    weighted GG &$6$  & \cite{priyadarshaniWGG}& 2015\\
    $\kappa-$GG &$4$ & \cite{2018kappaG,2021kappaGG} & 2018,2021\\
    Marshall-Olkin GG &$4$ & \cite{2018MOGG}& 2018\\
    double truncated GG &$5$ & \cite{2021DTGG} & 2021\\
     \hline
\end{tabularx}
    \caption{A timeline of the GG, its generalisations, and their number of parameters (including a scaling parameter).}
\label{tab:timeline}
\end{table}
Building on the success of the GG many authors have improved the applicability of the GG by the addition of parameters through generalisation. Generally speaking, the main excitement of these generalisations is their focus on superior fit in niche applications. In Table \ref{tab:timeline} a timeline of GG generalisations is given for completeness. The demands of new models today have a wider focus than simply better fit. The following authors \cite{ley14flex,jones15,mcleish1982robust,punzo2021multivariate,BTGN,ley2021flexible} specify these sometimes overlooked desirable qualities for new generalisations:
\begin{itemize}
     \item A low number of interpretable parameters: These include parameters that specifically control distribution shape qualities such as or similar to location, scale, skewness, and kurtosis.
     \item Favourable  estimation properties: It is important that the parameters can be estimated properly to ensure correct predictions and inferences from the model. Inferentially speaking, a generalised model for use in diagnosing distribution departure from a common baseline distribution.
     \item Simple mathematical tractability: Closed-form expressions and simple formulae aid in implementation, computational speed, and providing insight into the data.
    \item Finite moments: Most real-world measurements require this property.
\end{itemize}
Upon review of the literature in Table \ref{tab:timeline}, none of the generalisations except one, the $\kappa$-GG, have got parameters that are easily interpreted. This is due to these parameters having overlapping influence in distribution shape, which then obfuscates their role in achieving a certain fit. The $\kappa$-GG is such a distribution because it has an extra right-tail parameter which give geometric instead of exponential tails. Another drawback of these distributions is that they lack simple formulae. This is because of their complex setup of their generating mechanisms. The latter two points are seen as an opportunity for the derivation of a further generalisation of the GG distribution.

Here, we systematically construct a generalisation of the GG distribution that possesses interpretable parameters, favourable estimation, simple formulae, and finite moments because of its generating mechanism and setup. This is done by Weibullising the body-tail generalised normal distribution (BTN) \cite{BTGN} in order to have specific parameters for left-tail shape, body shape, and right-tail shape.

The paper is structured as follows. Section \ref{sec:genweibull} illustrates Weibullisation as a means to generate positive real line distributions from symmetric distributions. Section \ref{sec:FIG} introduces the flexible interpretable gamma distribution (FIG) generated by the Weibullisation of the BTN baseline distribution; The section further provides the derivations of the PDF, cumulative probability function (CDF), moments, moment generating function (MGF). Section \ref{sec:identifia} the parameters are proven to be identifiable. Section \ref{sec:estim} gives background on maximum likelihood (ML) estimation. Section \ref{sec:application} applies the FIG to hand grip strength and insurance loss data. Section \ref{sec:con} summarises the results and key findings.
\section*{Origins of gamma-like distributions}
In this section we investigate three different situations that give rise to gamma-like distributions. This insight may then be used to a guide an expansion of GG, with the aim of maintaining the shape clear and interpretable roles of the shape parameters.
\subsection*{Weibullisation}
\label{sec:genweibull}
\begin{figure}[h!]
  \includegraphics[width=13cm,keepaspectratio]{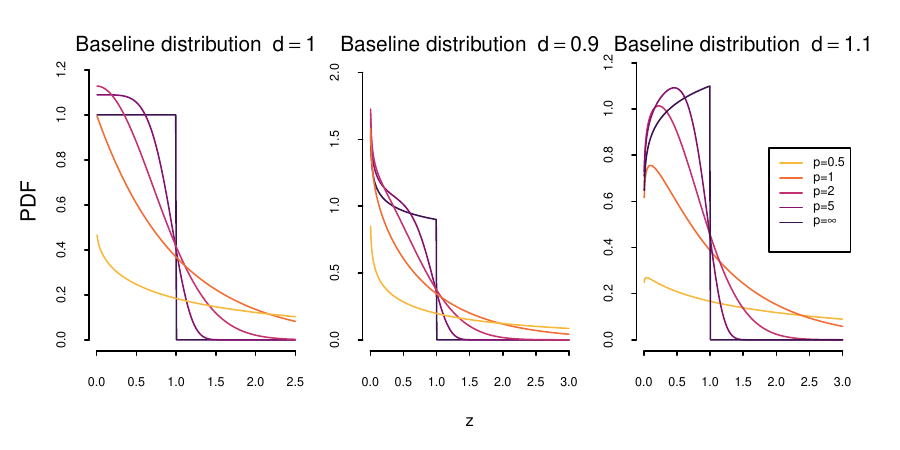}
 \caption{Examples of different Weibullisations of the half-GN equivalent to GG with different left-tail shapes $d$ and body shapes $p$.}
 \label{fig:gn_gg}
\end{figure}

The Weibullisation for a given baseline distribution of a random variable $Z$ occurs when considering the random variable $Z^{1/\nu}$ for $\nu>0$ \cite{PHAMGIA1989}. If the baseline distribution is selected to be symmetrical, the Weibullisation of the random variable $|Z|$ yields a positive distribution which has direct relationship to the shape of the baseline distribution. This process is illustrated for the GG distribution and its PDF in (\ref{eq:ggden}). As $|Z|$ is mathematically equivalent to one side or half of the baseline distribution due to symmetry, our analysis will continue with the latter. The GG distribution emerges when $Z$ follows a generalised normal distribution (GN); see \cite{subbotin23}. The PDF of the half-GN is given as:
\begin{equation}
 \label{eq:gn_den}
 f(z;s) = \frac{s}{\Gamma\left(\frac{1}{s}\right)}\text{e}^{-z^{s}},
\end{equation}
where $z,s>0$. Note that the half-GN contains the half-normal, half-Laplace, and uniform distributions for shape parameter values equal to $s=2$, $s=1$, and $s=\infty$, respectively \cite{subbotin23}. Let $y=z^{1/\nu}$, then the PDF of the random variable $Y=Z^{1/\nu}$ is given by:
\begin{equation}
 \label{eq:gngg}
 f(y;s,\nu) = \frac{s}{\Gamma\left(\frac{1}{s}\right)}y^{\nu-1}\text{e}^{-y^{s\nu}},
\end{equation}
where $y,s>0$. Comparing (\ref{eq:ggden}) and (\ref{eq:gngg}) shows that $Y$ follows a GG distribution with parameters $d=\nu$ and $p=s\nu$. In the kernel of (\ref{eq:gngg}), the right-tail is determined by $\text{e}^{-y^{s\nu}}$ and the left-tail by $y^{\nu-1}$, refer to Section \ref{sec:intro}. A visual representation of this generating process of (\ref{eq:gngg}) is given in Figure \ref{fig:gn_gg}. Observe that $d$ affects the left-tail shape, for $d=1$ we have no change to the half baseline distribution, for $d<1$ the left tail density is increased, and for $d>1$ the left-tail density is decreased. The left-tail shape does not influence the right tail, apart from a change of overall scale in the Weibullised distribution. The behaviour of the right tail is left to be determined by the half baseline distribution kernel from (\ref{eq:gn_den}) in this instance.
\subsection*{Power weighted distributions}
The GG distribution (\ref{eq:ggden}) also arises in the power weighted kernel of the integrand of the $r$-th absolute moments of the GN distribution. In general, if the absolute $r$-th moments of a distribution $Z$ exist and are finite, a positive distribution $Y$ can be generated from it. Consider the absolute $r$-th moment of $Z$:
\begin{equation}
\label{eq:absmom}
    E(|Z|^r) = \int_{\mathbb{R}}^{} |z|^{r} f(z) dz,
\end{equation}
where $r>0$. The integrand is a valid kernel for a positive support distribution, since its integral is finite by definition. Therefore, a new PDF can be generated by normalising the integrand function with the actual value of the integral, with the new PDF given by:
\begin{equation}
\label{eq:weibullshortcut}
    f(y;r) = \frac{y^{r} f(y)}{E(|Z|^{r})},
\end{equation}
where $y,r>0$ and $f(\cdot)$ is the original PDF of $Z$. To generate the GG PDF in this manner, $Z$ is taken as the GN distribution. By substituting (\ref{eq:gn_den}) into (\ref{eq:weibullshortcut}) the generated PDF of $Y$ is given by: 
\begin{equation}
\label{eq:ggdenmom}
     f(y;r,s) = \frac{s}{\Gamma\left(\frac{r}{s}\right)}y^{r}\text{e}^{-y^{s}},
\end{equation}
where $y,r,s>0$. Subsequently, from (\ref{eq:ggdenmom}), we have that $Y$ follows a GG distribution with $d=r+1$ and $p=s$.
\subsection*{Scale mixture model}
The FIG distribution also arises through a scale mixture of  power-function distribution (PF), which is a special case of the beta distribution \cite{johnson1970continuous}. Even this kind of mixture can not represent the GG, we still include it in this section, since it gives rise to the FIG distribution, which has gamma-like properties. The PDF of the scaled PF distribution is given below:
\begin{gather}
 f(x;u,\nu)= \frac{\nu}{u^{\nu}}x^{\nu-1},
\end{gather}
where $ 0 < x \leq u, u>0,\nu>0$, and is denoted as $X\sim PF(u,\nu)$.
\subsubsection*{Theorem:}
\label{thm:figstochastic}
Let $Z\sim PF(u,\nu)$ and $u\sim GG(\beta,\alpha)$, then the PDF of $Z$ is given by:
$f(z;\alpha,\beta,\nu)=\frac{\nu z^{\nu-1}}{\Gamma\left(\frac{\alpha+\nu}{\beta}\right)}\Gamma\left(\frac{\alpha}{\beta},z^\beta\right)$.
\subsubsection*{Proof}
Using the pre-defined random variables $Z$ and $U$. Noting that $Z<U$ by definition, and employing the indicator variable $I(\cdot)$ we have the following:
\begin{align}
    f(z;\alpha,\beta,\nu)&= \int^{}_{\mathbb{R}} \textit{I}\left(z\leq u\right)\frac{\nu}{ u^\nu}z^{\nu-1}\frac{u^{\alpha+\nu-1 }e^{-u^{\beta}}}{\Gamma\left(\frac{\alpha+\nu}{\beta}\right)} du\\&
    =\frac{\nu z^{\nu-1}}{\Gamma\left(\frac{\alpha+\nu}{\beta}\right)}\int^{}_{z}u^{\alpha-1}e^{-u^\beta}du\\&
    =\frac{\nu z^{\nu-1}}{\Gamma\left(\frac{\alpha+\nu}{\beta}\right)}\Gamma\left(\frac{\alpha}{\beta},z^\beta\right),
\end{align}
which concludes the proof.
\subsection*{Conclusion}
The selection of a baseline distribution, it's important to consider the existing roles of shape parameters and symmetry. This is crucial since any ambiguous roles pertaining to the baseline distribution parameters will inevitably be transferred to the distributions that are subsequently generated. Hence, it is not recommended to use asymmetric distributions due to the uncertainty of the effects of the asymmetry parameter following generalisation. This section presented three different methods of obtaining a GG type distributions. To derive the FIG distribution, our preference lies with the power weighted and scale mixture origins' parameterisations, details of which will be extensively discussed in the subsequent section.
\section*{The flexible interpretable gamma distribution}
\label{sec:FIG}
This section consists of the motivation for the chosen FIG baseline distribution and the derivations of the PDF, CDF, moments, and MGF for the standard and scaled FIG distribution.
\subsection*{Baseline distribution}
\begin{figure}[h!]
\centering
  \includegraphics[width=13cm,keepaspectratio]{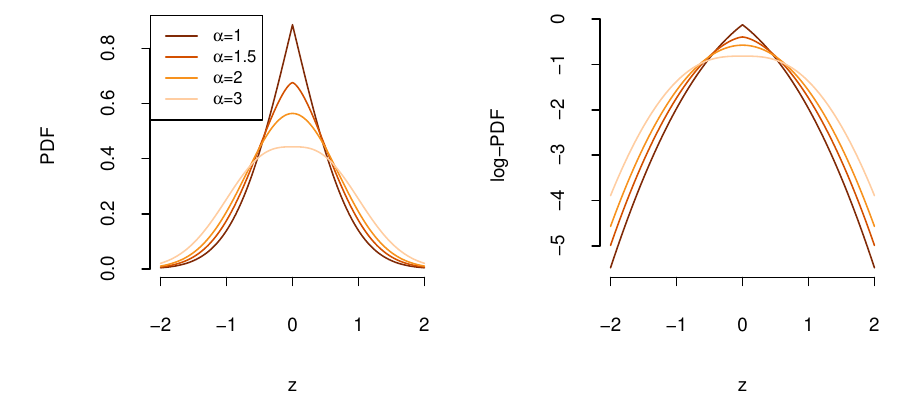}
 \caption{Examples of BTN PDFs for different values of body shape $\alpha$ and fixed tail  shape  $\beta=2$.}
 \label{fig:btgn_den_alpha}
\end{figure}

\begin{figure}[h!]
\centering
  \includegraphics[width=13cm,keepaspectratio]{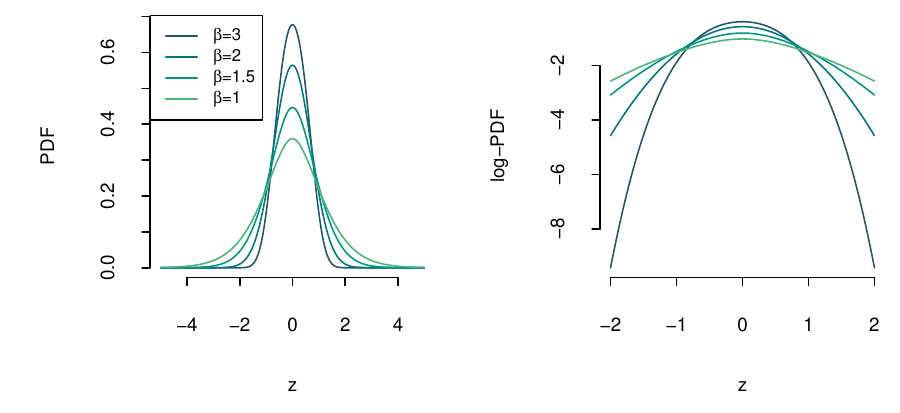}
 \caption{The BTN PDF for different values of tail shape $\beta$ and fixed body shape $\alpha=2$.}
 \label{fig:btgn_den_beta}
\end{figure}

The baseline distribution for the FIG is the BTN distribution. The BTN distribution is a generalisation of the GN and normal distribution which has interpretable parameters, simple mathematical tractability, and finite moments. The latter desirable properties will be transferred to the FIG distribution in the same way the Weibullisation of the GN transferred its properties to the GG. The PDF of the BTN is:
\begin{equation}
\label{eq:BTNpdf}
 f(z;\alpha,\beta)=\frac{\Gamma\left(\frac{\alpha}{\beta},|z|^\beta \right)}{2\Gamma\left( \frac{\alpha+1}{\beta}\right)},
\end{equation}
where $z\in\mathbb{R},\alpha, \beta>0$, and $\Gamma\left(\cdot,\cdot \right)$ is the upper incomplete gamma function; see (\cite{Gr<3dstheyn}, p. 899). The parameters have clear roles, where $\alpha$ determines body shape and $\beta$ determines the tail shape of the distribution. Note that, for $\alpha=\beta=s$, (\ref{eq:BTNpdf}) is equivalent to (\ref{eq:gn_den}), making the GN and its nested models a subset of the BTN; for more details refer to \cite{BTGN}. Due to the latter fact, the Weibullisation of the BTN will therefore contain the GG distribution for $\alpha=\beta$ as discussed in Section \ref{sec:genweibull}. The absolute moments of the BTN are given by: 
\begin{equation}
\label{eq:BTNmom}
E(|Z|^r) = \frac
 {\Gamma\left(\frac{\alpha+r+1}{\beta}\right)}
 {(r+1)\Gamma\left(\frac{\alpha+1}{\beta}\right)},
\end{equation}
where $r>0$; see \cite{BTGN}. In Figure \ref{fig:btgn_den_alpha} the different body shapes for a fixed tail shape can be seen. Similarly, in Figure \ref{fig:btgn_den_beta} the different tail shapes for a given body shape is shown.
The additional body shape parameter of the BTN specifically enhances the body shape of the GG distribution through Weibullisation in the FIG distribution. Therefore, the additional $\alpha$ parameter has an interpretation and provides information about the body shape of the FIG. Figure \ref{fig:FIG_GEN} illustrates this process of Weibullisation and the effect of the body shape parameter $\alpha$. Here, the GG and its fixed body shape is shown for $\alpha=\beta=2$. Notice that in the region of the body, $0.5<z<1$, the shape is determined by $\alpha$. In the region of the left tail, $z<0.5$, the shape is determined by $\nu$. In the region of the right tail, $z>1$, the shape is determined by $\beta$. Importantly, note that both the tail shapes stay markedly the same for different body shapes $\alpha$.
The FIG has extended the GG body shapes to ``steeper'' and ``flatter'' for body shape $\alpha<\beta$ and $\alpha>\beta$ while preserving the role of the left and right-tail parameters $\nu$ and $\beta$. In summary, the FIG distribution derives its unique properties from the BTN. The main property of interest is that each shape parameter specifically controls either left, body, or right-tail behaviour of the FIG. The FIG provides for greater flexibility in conjunction with numerically interpretable values for the cause of deviation from the GG.
\subsection*{PDF}
\begin{figure}[h!]
  \includegraphics[width=12cm,keepaspectratio]{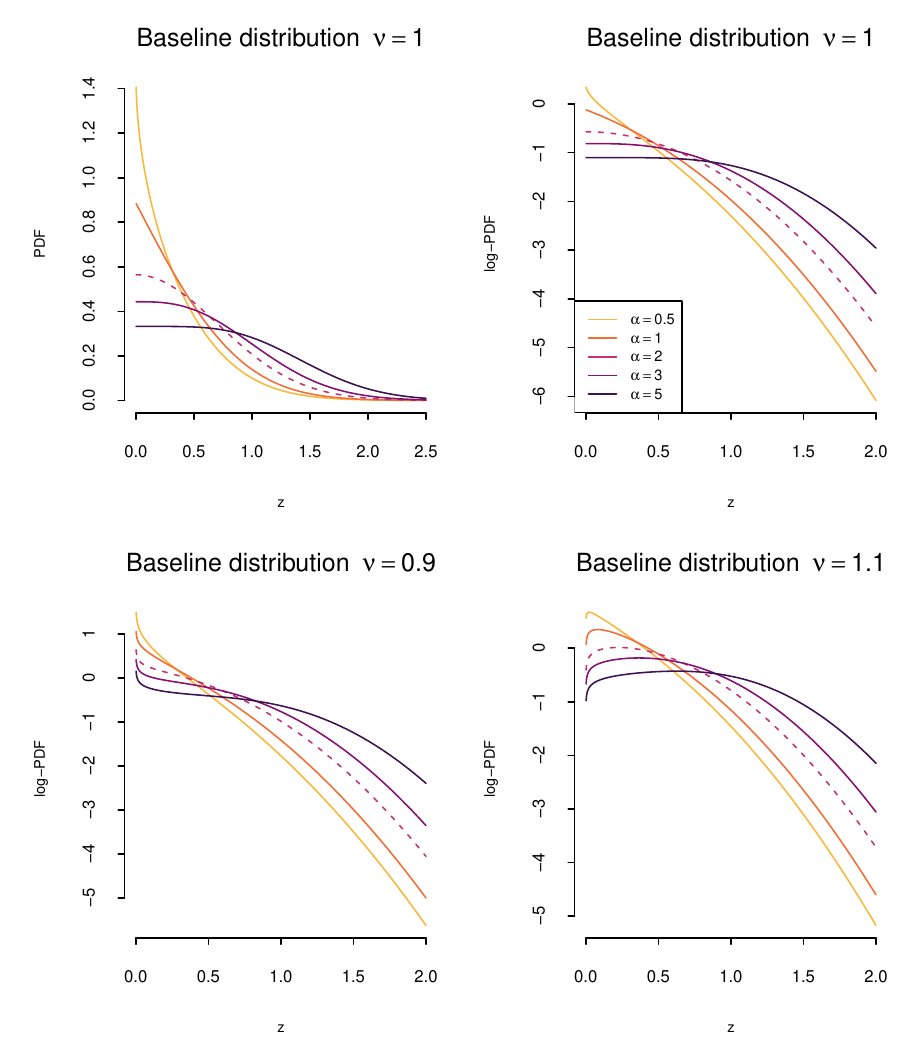}
 \caption{The baseline BTN and generated FIG PDFs for different body and left-tail shape parameters with fixed right-tail shape are shown in the first row (half baseline distributions) and second row (corresponding generated distributions).}
 \label{fig:FIG_GEN}
\end{figure}

The FIG PDF is derived by substituting (\ref{eq:BTNpdf}) and (\ref{eq:BTNmom}) into (\ref{eq:weibullshortcut}) and is given below:
\begin{equation}
\label{eq:FIGden}
    f(z;\alpha, \beta, \nu) =  
    \frac{\nu z^{\nu - 1}}{\Gamma \left( \frac{\alpha + \nu}{\beta} \right)}
  \Gamma \left( \frac{\alpha}{\beta}, z^{\beta} \right),
\end{equation}
where $z,\alpha, \beta, \nu>0$ denoted as $Z\sim FIG(\alpha,\beta,\nu)$. We note that if $\alpha=\beta$, then we have the standard GG distribution with parameters $p=\nu$ and $d=\alpha$. A depiction of the FIG PDF (\ref{eq:FIGden}) and corresponding baseline PDF is given in Figure \ref{fig:FIG_GEN}. The PDF of the scaled FIG,  denoted as $X\sim FIG(\sigma,\alpha,\beta,\nu)$, is obtained using the transformation $X = \sigma Z$:
\begin{align}
 f(x; \sigma, \alpha, \beta, \nu ) & =  
 \frac{\nu x^{\nu - 1}}{\sigma^\nu
 \Gamma \left( \frac{\alpha + \nu}{\beta} \right)}
  \Gamma \left( \frac{\alpha}{\beta}, \left(\frac{x}{\sigma}\right)^{\beta} \right),
  \label{eq:dendenscaleFIG}
\end{align}
where $\sigma>0$. Since the scaled FIG is a generalisation of the GG distribution it has many sub-models which is summarised in Table \ref{tab:submodels}.
\begin{center}
\begin{table}[h!]
    \centering
    \begin{tabular}{| c c c c c |}
    \hline
    Distribution & $\sigma$ & $\alpha$ & $\beta$ & $\nu$ \\ 
    \hline
    chi-squared & 2 & 1 & 1 & $df/2$ \\
    exponential  &$1/\lambda$ & 1 & 1 & 1\\
    gamma       & $\theta$ &  1  &  1  & $k$\\
    generalised gamma & $a$ & $p$ & $p$ & $d$\\
    half-body-tail generalised normal &$\sigma$ & $\alpha$ & $\beta$ & 1 \\
    half-generalised normal &$\sigma$ & s & s & 1 \\
    half-normal  &$\sqrt{2}\sigma$ & 2 & 2 &  1\\
    Maxwell–Boltzmann & $\sqrt{2}a$ & 2 & 2 & 3 \\
    Rayleigh & $\sqrt{2}\sigma$ & 2 & 2 & 2 \\
    uniform & $\cdot$ & $\infty$ & $\infty$ & 1 \\
    Weibull & $\lambda$ & $k$ & $k$ & $k$\\
    \hline
    \end{tabular}
    \caption{A summary of nested models of the FIG distribution.}
\label{tab:submodels}
\end{table}
\end{center}
The CDF of the standard FIG is derived with the definition of a CDF and (\ref{eq:FIGden}):
\begin{align}
 F(z; \alpha, \beta, \nu) & = \int_{-\infty}^{z} 
 \frac{\nu t^{\nu - 1}}{\Gamma \left( \frac{\alpha + \nu}{\beta} \right)}
  \Gamma \left( \frac{\alpha}{\beta}, t^{\beta} \right) dt
 \notag\\ 
 &=
 \frac{\nu}{\Gamma \left( \frac{\alpha + \nu}{\beta} \right)}
 \left(
 \int_{0}^{\infty} t^{\nu - 1} \Gamma \left( \frac{\alpha}{\beta}, t^{\beta} \right) dt
 - \int_{z}^{\infty} t^{\nu - 1} \Gamma \left( \frac{\alpha}{\beta}, t^{\beta} \right) dt
 \right).\notag
\end{align}
Applying Lemma 2 of the Appendix to both integrals, we have that:
\begin{align}
  F(z; \alpha, \beta, \nu)  & =
 \frac{\nu}{\Gamma \left( \frac{\alpha + \nu}{\beta} \right)}
\left(
 \lim_{t \to 0^+} \frac{
  \Gamma \left( \frac{\alpha + \nu}{\beta}, t^{\beta} \right)
t^{\nu} \Gamma \left( \frac{\alpha}{\beta}, t^{\beta} \right)
 }{\nu}
 - \frac{
  \Gamma \left( \frac{\alpha + \nu}{\beta}, z^{\beta} \right)
  - z^{\nu} \Gamma \left( \frac{\alpha}{\beta}, z^{\beta} \right)
 }{\nu}
 \right)\notag\\
 &=
 \frac{
  \gamma \left( \frac{\alpha + \nu}{\beta}, z^{\beta} \right)
  + z^{\nu} \Gamma \left( \frac{\alpha}{\beta}, z^{\beta} \right)
 }{\Gamma \left( \frac{\alpha + \nu}{\beta} \right)},\label{eq:FIGcdf}
\end{align}
where $\gamma(\cdot,\cdot)$ is the lower incomplete gamma function; see (\cite{Gr<3dstheyn}, p. 899). Subsequently, the CDF of $X\sim FIG(\sigma,\alpha,\beta,\nu)$ is given by the substitution of $z=\frac{x}{\sigma}$ in (\ref{eq:FIGcdf}).
\subsection*{Mode}
\label{thm:unimode}
The maximum of the standard {FIG} {PDF} is given by the maximum of the {FIG} kernel in (\ref{eq:FIGden}). We consider two cases for obtaining the mode of the FIG. For $\nu\leq 1$ we have that
\begin{equation}
    \lim_{z \to 0^+} z^{\nu-1}\Gamma\left(\frac{\alpha}{\beta},z^{\beta}\right)=\infty,
\end{equation}
and
\begin{equation}
    \lim_{z \to \infty} z^{\nu-1}\Gamma\left(\frac{\alpha}{\beta},z^{\beta}\right)=0,
\end{equation}
which implies that the mode is zero. For $\nu>1$ we have that
\begin{align}
    \frac{\partial}{\partial z} z^{\nu-1}\Gamma\left(\frac{\alpha}{\beta},z^{\beta}\right)&=
    (\nu-1)z^{\nu-2}\Gamma\left(\frac{\alpha}{\beta},z^{\beta}\right)
    +z^{\nu-1}\left(-z^{\beta(\frac{\alpha}{\beta}-1)}e^{-z^{\beta}}\beta z^{\beta-1}\right)\notag\\
    &= (\nu-1)z^{\nu-2}\Gamma\left(\frac{\alpha}{\beta},z^{\beta}\right)
    -z^{\nu-2}\left(z^{\alpha}e^{-z^{\beta}}\beta \right)\notag\\
    &=z^{\nu-2}
    \left((\nu-1)\Gamma\left(\frac{\alpha}{\beta},z^{\beta}\right)
    -\beta z^{\alpha}e^{-z^{\beta}}\right) \label{eq:modemax2}.
\end{align}
Examining the elements inside the brackets in (\ref{eq:modemax2}), we find that for
\begin{equation}
    \lim_{z\to 0^+}(\nu-1)\Gamma\left(\frac{\alpha}{\beta},z^{\beta}\right)=(\nu-1)\Gamma\left(\frac{\alpha}{\beta}\right)
    >
    \lim_{z\to 0^+}\beta z^{\alpha}e^{-z^{\beta}}=0, \label{eq:modegthan}
\end{equation}
which implies (\ref{eq:modemax2}) is greater than zero for some $z$. Sine both terms in (\ref{eq:modemax2}) are greater than zero for all $z$, we investigate the rate of decrease for these function as $z$ increases. The relative rate of decrease for the left and right terms in (\ref{eq:modemax2}) with respect to $z$ is given by
\begin{equation}
\frac
{\frac{\partial}{\partial z}\beta z^{\alpha}e^{-z^{\beta}}}
{\frac{\partial}{\partial z}(\nu-1)\Gamma\left(\frac{\alpha}{\beta},z^{\beta}\right)}
=
\frac
{\beta \left(\alpha z^{\alpha-1}e^{-z^{\beta}}+z^{\alpha}e^{-z^{\beta}}\left(-\beta z^{\beta-1}\right)\right)}
{(\nu-1)\left(-\beta z^{\alpha-1}e^{-z^{\beta}}\right)}
=
\frac
{\beta z^{\beta}-\alpha}
{\nu-1},
\end{equation}
which implies that $\beta z^{\alpha}e^{-z^{\beta}}$ decrease geometrically faster than $(\nu-1)\Gamma\left(\frac{\alpha}{\beta},z^{\beta}\right)$ for $z>\frac{\text{ln}(\alpha)-\text{ln}(\beta)}{\beta}$. Noting that $(\nu-1)\Gamma\left(\frac{\alpha}{\beta},z^{\beta}\right)$ is monotonically decreasing, and $\beta z^{\alpha}e^{-z^{\beta}}$ increasing and then decreasing but at a slower relative rate for all $z$ we have that (\ref{eq:modemax2}) is negative from some $z$ onward. Therefore, we have that PDF of the is first increasing and then decreasing for $\nu>1$ implying a mode which can be numerically calculated by setting (\ref{eq:modemax2}) to zero.
\subsection*{Tail behaviour}
To examine the behaviour of the left and right tail of the FIG in comparison to the GG, we examine the derivative of the difference of log-kernel functions between the FIG and GG as follows:
\begin{align}
    d(z;\alpha,\beta,\nu)&=\frac{\partial}{\partial z}\left(\ln{\left(z^{\nu-1}\Gamma\left(\frac{\alpha}{\beta}, z^\beta\right)\right)} -\ln{\left(z^{\nu-1}e^{-z^\beta}\right)}\right)\notag\\&
    =\frac{\partial}{\partial z}\left(\ln{\left(\Gamma\left(\frac{\alpha}{\beta}, z^\beta\right)\right)} -{z^\beta}\right)\notag\\&
    =\frac{-\left(z^\beta\right)^{\frac{\alpha}{\beta}-1}e^{-z^\beta}\left(-\beta z^{\beta-1}\right)}{\Gamma\left(\frac{\alpha}{\beta}, z^\beta\right)}-\beta z^{\beta-1}\notag\\&
    =\frac{-\beta z^{\alpha-1}e^{-z^\beta}}{\Gamma\left(\frac{\alpha}{\beta}, z^\beta\right)}-\beta z^{\beta-1}. \label{eq:figlogderdiff}
\end{align}
Comparing the left tail behaviour of the FIG to the GG, we evaluate $\lim_{z \to 0^+}d(z;\alpha,\beta,\nu)={0\times1}-0 = 0$, which suggests that the left tail behaviour of the FIG approximates the shape of the GG distribution for small $z$. This property is visually confirmed by the bottom left sub-figure of Figure \ref{fig:FIG_GEN}. Regarding, the right tail behaviour, we first concentrate on the first term in (\ref{eq:figlogderdiff}) for large values of $z$. The limit is of the form zero divided by zero, for which we apply L`Hospitals rule:
\begin{align}
    &\lim_{z \to \infty}\frac{-\beta z^{\alpha-1}e^{-z^\beta}}{\Gamma\left(\frac{\alpha}{\beta}, z^\beta\right)}\notag\\
    &\stackrel{L'H}{=}\lim_{z \to \infty}\frac{-\beta\left((\alpha-1)z^{\alpha-z}e^{-z^\beta}+z^{\alpha-1}e^{-z^\beta}\left(-\beta z^{\beta-1}\right)\right)}{-\beta z^{\alpha-1}e^{-z^{\beta}}}\notag\\
    &=\lim_{z \to \infty}(\alpha-1)z^{-1}-z^{\beta-1}\notag\\
    &=-z^{\beta-1}.\label{eq:righttailbehavior}
\end{align}
Substituting the result from (\ref{eq:righttailbehavior}) into $\lim_{z \to \infty}d(z;\alpha,\beta,\nu)=0$, we conclude that the FIG right tail behaviour approximates the shape of the GG distribution for large $z$. Similarly, this property is visually confirmed by the bottom right sub-figure of Figure \ref{fig:FIG_GEN}. It is worth noting that the body shape parameter $\alpha$ does not influence calculated left and right tail limits, suggesting that the left and right tail parameters have maintained their roles and interpretation.
\subsection*{Moments}
The $r$-th moment of the standard FIG is derived from (\ref{eq:FIGden}), and Lemma 3 in the Appendix:
\begin{equation}
E(Z^r)=
\frac{\nu \Gamma \left( \frac{\alpha + \nu + r}{\beta} \right)}{(\nu + r) \Gamma \left( \frac{\alpha + \nu}{\beta} \right)}.\label{eq:rawmom}
\end{equation}
Subsequently, the $r$th moment of $X\sim FIG(\sigma,\alpha,\beta,\nu)$ is given by (\ref{eq:rawmom}), and the identity $E(X^r)=\sigma^r E(Z^r)$.
\subsection*{Moment generating function}
From the definition of a MGF, (\ref{eq:FIGden}) and the series expansion of the incomplete gamma function (\cite{Gr<3dstheyn}, p. 901). The MGF of the standard FIG is derived by:
\begin{equation}
M_Z(t)  = \frac{\nu}{t^\nu \Gamma\left(\frac{\alpha+\nu}{\beta}\right)}
\left(
\Gamma \left(\frac{\alpha}{\beta}\right)
\Gamma \left(\nu\right)
-
\sum_{k=0}^{\infty}
            \frac{(-1)^k \Gamma \left(\alpha+\beta k+\nu\right)}
            {k!(\alpha/\beta+k)t^{\alpha+\beta k}}
            dz\right).
            \label{eq:mgf}
\end{equation}
Subsequently, the MGF of $X\sim BTN(\mu,\sigma,\alpha,\beta)$ is given by (\ref{eq:mgf}), and the identity $M_X(t) =\text{e}^{t\mu}M_Z(t\sigma)$.
\section*{FIG identifiability}
\label{sec:identifia}
The FIG distribution is derived with the intent of gaining insights from the fitted parameters to data. It is, therefore, important that the parameters of the FIG PDF (\ref{eq:dendenscaleFIG}) are mathematically identifiable.
\subsection*{Theorem:}
\label{thm:identi}
Let $X\sim FIG(\sigma, \alpha, \beta,  \nu)$. If $\sigma_{1}, \sigma_{2}>0, \alpha_{1}, \alpha_{2}>0, \beta_{1}, \beta_{2}>0,  \nu_{1},  \nu_{2}>0$ such that $f\left(x ; \sigma_{1}, \alpha_{1}, \beta_{1},  \nu_{1}\right)=f\left(x ; \sigma_{2}, \alpha_{2}, \beta_{2},  \nu_{2}\right), \forall x >0$ then $\sigma_{1}=\sigma_{2}, \alpha_{1}=\alpha_{2}, \beta_{1}=\beta_{2},  \nu_{1}= \nu_{2}$.
\subsection*{Proof}
The proof follows by method of contradiction, assume the parameters of the FIN are not identifiable. That is, there exist parameters $\sigma_1\neq\sigma_2,\alpha_1\neq\alpha_2,\beta_1\neq\beta_2,\nu_1\neq\nu_2$ such that 
\begin{equation}
\label{eq:identhyp}
f\left(x ;\sigma_{1}, \alpha_{1}, \beta_{1},  \nu_{1}\right)=f\left(x ;\sigma_{2}, \alpha_{2}, \beta_{2},  \nu_{2}\right),
\end{equation}
where $\sigma_1,\sigma_2,\alpha_1,\alpha_2,\beta_1,\beta_2,\nu_1,\nu_2>0$. From the hypothesis (\ref{eq:identhyp}) and substitution of (\ref{eq:dendenscaleFIG}) we have that
\begin{align}
&\frac{ \nu_{1} x^{ \nu_{1}-1} \Gamma\left(\frac{\alpha_{1}}{\beta_{1}},\left(\frac{x}{\sigma_{1}}\right)^{\beta_{1}}\right)}{\sigma_{1}^{ \nu_{1}} \Gamma\left(\frac{\alpha_{1}+ \nu_{1}}{\beta_{1}}\right)}=\frac{ \nu_{2} x^{ \nu_{2}-1} \Gamma\left(\frac{\alpha_{2}}{\beta_{2}},\left(\frac{x}{\sigma_{2}}\right)^{\beta_{2}}\right)}{\sigma_{2}^{ \nu_{2}} \Gamma\left(\frac{\alpha_{2}+ \nu_{2}}{\beta_{2}}\right)}, \forall x >0,\notag
\end{align}
from which we obtain
\begin{equation}
\frac{x^{ \nu_{1}- \nu_{2}} \Gamma\left(\frac{\alpha_{1}}{\beta_{1}},\left(\frac{x}{\sigma_{1}}\right)^{\beta_{1}}\right)}{\Gamma\left(\frac{\alpha_{2}}{\beta_{2}},\left(\frac{x}{\sigma_{2}}\right)^{\beta_{2}}\right)}=\frac{ \nu_{2} \sigma_{1}^{ \nu_{1}} \Gamma\left(\frac{\alpha_{1}+ \nu_{1}}{\beta_{1}}\right)}{ \nu_{1} \sigma_{2}  \nu_{2} \Gamma\left(\frac{\alpha_{2}+ \nu_{2}}{\beta_{2}}\right)},  \forall x >0.
\label{eq:nucontra}
\end{equation}
Assuming for contradiction that $ \nu_{1} \neq  \nu_{2}$. It would now be possible that
\begin{equation}
\lim _{x \rightarrow 0^{+}} \frac{x^{ \nu_{1}- \nu_{2}} \Gamma\left(\frac{\alpha_{1}}{\beta_{1}},\left(\frac{x}{\sigma_{1}}\right)^{\beta_{1}}\right)}{\Gamma\left(\frac{\alpha_{2}}{\beta_{2}},\left(\frac{x}{\sigma_{2}}\right)^{\beta_{2}}\right)}=\left\{\begin{array}{lll}
0^{+} & \text {se } &  \nu_{1}> \nu_{2} \\
+\infty & \text { se } &  \nu_{1}< \nu_{2},
\end{array}\right.\notag
\end{equation}
which is in contradiction with the equality in (\ref{eq:identhyp}) where $\frac{ \nu_{2} \sigma_{1}^{ \nu_{1}} \Gamma\left(\frac{\alpha_{1}+ \nu_{1}}{\beta_{1}}\right)}{ \nu_{1} \sigma_{2}  \nu_{2} \Gamma\left(\frac{\alpha_{2}+ \nu_{2}}{\beta_{2}}\right)} \in \mathbb{R}^{+}$. Therefore, it must necessarily be that ${ \nu}_{\mathbf{1}}={ \nu}_{\mathbf{2}}$. We therefore substitute $ \nu= \nu_{1}= \nu_{2}$ from here on forward. From (\ref{eq:identhyp}), it follows that:
\begin{gather}
\Gamma\left(\frac{\alpha_{1}}{\beta_{1}},\left(\frac{x}{\sigma_{1}}\right)^{\beta_{1}}\right)=\frac{\sigma_{1}^{ \nu} \Gamma\left(\frac{\alpha_{1}+ \nu}{\beta_{1}}\right)}{\sigma_{2}^{ \nu} \Gamma\left(\frac{\alpha_{2}+ \nu}{\beta_{2}}\right)} \Gamma\left(\frac{\alpha_{2}}{\beta_{2}},\left(\frac{x}{\sigma_{2}}\right)^{\beta_{2}}\right),  \forall x >0, \notag\\
\Gamma\left(\frac{\alpha_{1}}{\beta_{1}},\left(\frac{x}{\sigma_{1}}\right)^{\beta_{1}}\right)=k \Gamma\left(\frac{\alpha_{2}}{\beta_{2}},\left(\frac{x}{\sigma_{2}}\right)^{\beta_{2}}\right),  \forall x >0,
\label{eq:ksimplified}
\end{gather}
where
$
k=\frac{\sigma_{1}^{ \nu} \Gamma\left(\frac{\alpha_{1}+ \nu}{\beta_{1}}\right)}{\sigma_{2}^{ \nu} \Gamma\left(\frac{\alpha_{2}+ \nu}{\beta_{2}}\right)}>0$.
Taking the derivative of both sides of (\ref{eq:ksimplified}) with respect to $x$, we obtain
\begin{equation}
-\left(\frac{x}{\sigma_{1}}\right)^{\beta_{1}\left(\frac{\alpha_{1}}{\beta_{1}}-1\right)}\text{e}^{-\left(\frac{x}{\sigma_{1}}\right)^{\beta_{1}}} \frac{\beta_{1}}{\sigma_{1}}\left(\frac{x}{\sigma_{1}}\right)^{\beta_{1}-1}=-k\left(\frac{x}{\sigma_{2}}\right)^{\beta_{2}\left(\frac{\alpha_{2}}{\beta_{2}}-1\right)}\text{e}^{-\left(\frac{x}{\sigma_{2}}\right)^{\beta_{2}} \frac{\beta_{2}}{\sigma_{2}}\left(\frac{x}{\sigma_{2}}\right)^{\beta_{2}-1}},\forall x >0.\notag
\end{equation}
After rearranging, it follows that
\begin{gather}
x^{\alpha_{1}-\alpha_{2}}\text{e}^{\left(\frac{x}{\sigma_{2}}\right)^{\beta_{2}}-\left(\frac{x}{\sigma_{1}}\right)^{\beta_{1}}}=\frac{k \beta_{2} \sigma_{1}^{\alpha_{1}}}{\beta_{1} \sigma_{2}^{\alpha_{2}}},  \forall x >0.
\label{eq:contra2}
\end{gather}
Assuming for contradiction that $\alpha_{1} \neq \alpha_{2}$. It would now be possible that
$$
\lim _{x \rightarrow 0^{+}} x^{\alpha_{1}-\alpha_{2}}\text{e}^{\left(\frac{x}{\sigma_{2}}\right)^{\beta_{2}}-\left(\frac{x}{\sigma_{1}}\right)^{\beta_{1}}}=\left\{\begin{array}{lll}
0^{+} & \text {if} & \alpha_{1}>\alpha_{2} \\
+\infty & \text {if} & \alpha_{1}<\alpha_{2},
\end{array}\right.
$$
which is in contradiction with the equality in (\ref{eq:contra2}), since $\frac{k \beta_{2} \sigma_{1}{ }^{\alpha_{1}}}{\beta_{1} \sigma_{2} \alpha_{2}} \in \mathbb{R}^{+}$. Therefore, it must necessarily be that ${\alpha}_{\mathbf{1}}={\alpha}_{\mathbf{2}}$. We therefore substitute $\alpha=\alpha_{1}=\alpha_{2}$. Consequently, (\ref{eq:contra2}) simplifies to:
\begin{equation}
   \label{eq:contra3}
 \text{e}^{\left(\frac{x}{\sigma_{2}}\right)^{\beta_{2}}-\left(\frac{x}{\sigma_{1}}\right)^{\beta_{1}}}=\frac{k \beta_{2} \sigma_{1}^{\alpha}}{\beta_{1} \sigma_{2}^{\alpha}}, \forall x >0.
\end{equation}
Assuming for contradiction that $\beta_{1} \neq \beta_{2}$. It would now be possible that
$$
\lim _{x \rightarrow+\infty}\text{e}^{\left(\frac{x}{\sigma_{2}}\right)^{\beta_{2}}-\left(\frac{x}{\sigma_{1}}\right)^{\beta_{1}}}=\left\{\begin{array}{lll}
0^{+} & \text {if} & \beta_{1}>\beta_{2} \\
+\infty & \text {if} & \beta_{1}<\beta_{2},
\end{array}\right.
$$
which is in contradiction with the equality in (\ref{eq:contra3}), since $\frac{k \beta_{2} \sigma_{1}{ }^{\alpha}}{\beta_{1} \sigma_{2}{ }^{\alpha}} \in \mathbb{R}^{+}$. Therefore, it must necessarily be that ${\beta}_{1}={\beta}_{2}$. We therefore substitute $\beta=\beta_{1}=\beta_{2}$ from here on forward. Consequently, (\ref{eq:contra3}) simplifies to:
\begin{equation}
\label{eq:contra4}
 \text{e}^{\left(\frac{x}{\sigma_{2}}\right)^{\beta}-\left(\frac{x}{\sigma_{1}}\right)^{\beta}}=\frac{k \sigma_{1}^{\alpha}}{\sigma_{2}^{\alpha}}, \forall x >0.
\end{equation}
Assuming for contradiction that $\sigma_{1} \neq \sigma_{2}$. It would now be possible that
$$
\lim _{x \rightarrow+\infty}\text{e}^{\left(\frac{x}{\sigma_{2}}\right)^{\beta}-\left(\frac{x}{\sigma_{1}}\right)^{\beta}}=\lim _{x \rightarrow+\infty}\text{e}^{\frac{\sigma_{1}^{\beta}-\sigma_{2} \beta}{\sigma_{2} \beta \sigma_{1} \beta} x^{\beta}}=\left\{\begin{array}{lll}
+\infty & \text {if} & \sigma_{1}>\sigma_{2} \\
0^{+} & \text {if} & \sigma_{1}<\sigma_{2},
\end{array}\right.
$$
which is in contradiction with the equality in (\ref{eq:contra4}), since $\frac{k \sigma_{1}^{\alpha}}{\sigma_{2}^{\alpha}} \in \mathbb{R}^{+}$. Therefore, it must necessarily be that $\sigma_{1}=\sigma_{2}$. In summary, it has been proven that $\sigma_{1}=\sigma_{2}, \alpha_{1}=\alpha_{2}, \beta_{1}=\beta_{2}, \nu_{1}= \nu_{2}$ which completes the proof.

\section*{FIG maximum likelihood equations}
\label{sec:estim}
The log-likelihood (LL) for a random sample $x_1,x_2,\dots,x_n$ from $X\sim FIG(\sigma,\alpha,\beta,\nu)$ observations is
\begin{align}\label{eq:LL}
 \text{LL}(\sigma, \alpha, \beta, \nu;x_1,x_2,\dots,x_n) & =
 \sum_{i=1}^{n}
 \Bigg(
 \left.
 \ln{\left( \nu \right)} -\ln(\sigma) 
 + (\nu - 1)\ln{( z_i )} \right.\notag\\
 &\phantom{\sum} \left.
 + \ln{\left( \Gamma \left( \frac{\alpha}{\beta}, z^\beta \right) \right)} 
 - \ln{\left( \Gamma \left( \frac{\alpha + \nu}{\beta} \right) \right)}
 \right),
\end{align}
where $z_i=x_i/\sigma$. The derivatives of the individual terms in (\ref{eq:LL}) with respect to the FIG parameters are given by:
\begin{align*}
\frac{\partial LL}{\partial \sigma}
&=
-\frac{\nu}{\sigma}
-\frac{\beta z_i^{\beta}\rpsi_2\left(\frac{\alpha}{\beta},z_i^\beta\right)}
{\sigma\Gamma\left(\frac{\alpha}{\beta},z_i^\beta\right)}\\
 \frac{\partial LL}{\partial\alpha} & = 
 \frac{\rpsi_1 \left(\frac{\alpha}{\beta},z_i^{\beta}\right)}{\beta\Gamma \left( \frac{\alpha}{\beta}, z_i^\beta \right) }
 -\frac{\psi\left(\frac{\alpha+\nu}{\beta}\right)}{\beta}\\
 \frac{\partial LL}{\partial\beta}
 & =\frac{\rpsi_3 \left(\frac{\alpha}{\beta},z_i^{\beta},\beta\right)}{\Gamma \left( \frac{\alpha}{\beta}, z_i^\beta \right)}
 +\frac{\psi\left(\frac{\alpha+\nu}{\beta}\right)\left(\alpha+\nu\right)}{\beta^{2}}\\
 \frac{\partial LL}{\partial \nu}
 & = \frac{1}{\nu} + \ln{\left(z_i\right)} -\frac{\psi\left(\frac{\alpha+\nu}{\beta}\right)}{\beta},
\end{align*}
where
\begin{align*}
 \rpsi_1(u,v)  & = \frac{\partial}{\partial u} \Gamma(u,v) = \Gamma \left(u,v\right)\ln{v}+A(u,v),    \\
 \rpsi_2(u,v)  & = \frac{\partial}{\partial v} \Gamma(u,v) = -v^{u-1}\text{e}^{-v},            \\
 \rpsi_3(u,v,w) & = \rpsi_1(u,v)\frac{\partial}{\partial w} u+ \rpsi_2(u,v)\frac{\partial}{\partial w} v ,\\
 A(u, v)    & = G^{3,0}_{2,3} \left(v \left \lvert \begin{gathered} 1, 1 \\ 0, 0, u \end{gathered} \right. \right),
\end{align*}
$\psi(\cdot)$ is the digamma function, and $G$ is the Meijer's G function; see (\cite{Gr<3dstheyn}, p. 850,902).
\section*{FIG applications}
\label{sec:application}
In this section the FIG is applied to commonly available benchmark data to compare the flexibility with competing distributions. The competitor models are the inverse Guassian (IG) \cite{IGbook} and generalised inverse Guassian (GIG) \cite{GIGbook}. These models are some of the most famous for modelling positive data; see \cite{tweedie1957statistical} and review paper \cite{folks1978inverse}. For a detailed list  where the IG has been successfully implemented; see \cite{seshadri2012inverse} and \cite{johnson1970continuous}. The evaluation of fit is done by computing both in-sample and out-of-sample validation metrics. The in-sample statistics are Akaike information criterion ($AIC_{is}$) and Bayesian information criterion ($BIC_{is}$) computed on the subset of data used for estimation; see \cite{akaike1974new,schwarz1978estimating}. The out-of-sample LL ($LL_{os}$) is computed on a $10\%$ subset of data excluded from estimation. This is done to ensure robust goodness of fit analysis and the prevention of over fit of the final models. The application is implemented using packages \cite{gamlss}, NumPy \cite{numpy}, Scipy \cite{scipy}, and mpmath \cite{mpmath} in R and Python.
\subsection*{Hand grip strength}
The data consists of the hand grip strength English school boys available in the gammlss.data package, available online at \url{https://cran.r-project.org/web/packages/gamlss.data}, accessed 23 July 2022. The summary statistics of hand grip strength are given in Table \ref{tab:gripsummary}. 
\begin{table}[h!]
  \centering
  \begin{tabular}{cccccccc}
  	\hline
        {Min}
  	&{Max}
  	&{Median}
  	&{Mean}
  	&{{Std}}
  	&{\makecell{Pearson \\Skewness}}
  	&{\makecell{Pearson \\Kurtosis}}
  	&{Count}\\
  	\hline
  	7& 60& 24& 25.71 & 8.85 & 0.81 & 3.41 & 3766\\
  	\hline
  \end{tabular}
  \caption{Summary statistics for hand grip strength data.}
  \label{tab:gripsummary}
\end{table}
\begin{table}[h!]
  \centering
  \begin{tabular}{lrrrr}
  \hline
   Distribution &   $LL_{is}$  &    $AIC_{is}$ &    $BIC_{is}$ &   $LL_{os}$ \\
  \hline
  IG & -11974.224 & 23952.447 & 23964.704 & -1352.182 \\ 
  GIG & -11973.171 & 23952.342 & 23970.728 & -1349.809 \\ 
  GG & -11982.560 & 23971.120 & 23989.506 & -1342.772 \\
  FIG & -11976.559 & 23961.118 & 23985.633 & -1344.217 \\
  \hline
  \end{tabular}
  \caption{In- and out-of-sample metrics of distributions fitted to hand grip strength data.}
  \label{tab:gripmetrics}
\end{table}
The in-sample criterion and out-of-sample $LL_{os}$ are tabulated in Table \ref{tab:gripmetrics}. In this application, the in-sample metrics are the lowest for the IG and GIG. However, the $LL_{os}$ favours the GG and FIG because they are higher than the IG and GIG. It can therefore be concluded that the distributions perform similarly, with preference to be given to the simpler GG and IG due to parsimony. The fitting of the generalised models as competitors remain important, since we would not know whether a more complex model is necessary if we do not fit one.
\subsection*{Danish Fire losses}
\begin{figure}[h!]
\centering
 \includegraphics[width=10cm,keepaspectratio]{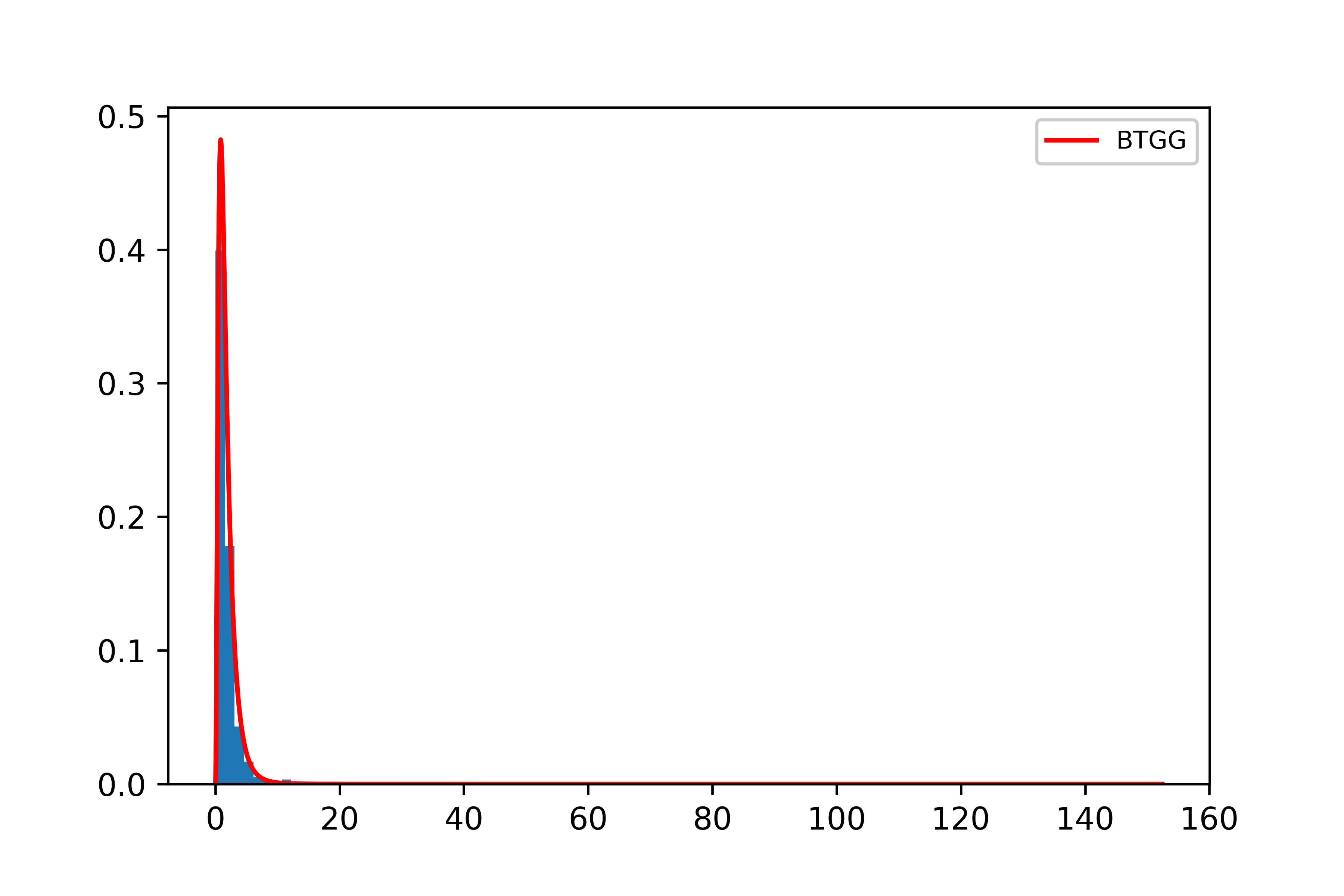}
 \caption{Relative frequency and fitted FIG PDF to Danish building fire loss data.}
 \label{fig:danishbuildingsfit}
\end{figure}

The data consists of the Copenhagen Reinsurance fire losses for the period from 1980 to 1990. The claim amount in millions of Danish Krone are divided into building, contents, and profit loss, available in the CASdatasets package, available online at \url{http://cas.uqam.ca/}, accessed 25 July 2022. The summary statistics of the different loss types are given in Table \ref{tab:danishsummary}. An illustration of the distributional shape shown by the relative frequency and fitted FIG PDF is shown in Figure \ref{fig:danishbuildingsfit}. 

\begin{table}[h!]
  \centering
  \begin{tabular}{ccrcccccr}
  	\hline
  	Type
    &{Min}
  	&{Max}
  	&{Median}
  	&{Mean}
  	&{{Std}}
  	&{\makecell{Pearson \\Skewness}}
  	&{\makecell{Pearson \\Kurtosis}}
  	&{Count}\\
  	\hline
    Buildings&0.02 & 152.41 & 1.33 & 2.01 & 4.72 & 22.98 & 660.4 & 1990\\
    Content&0.01 & 132.01 & 0.58 & 1.71 & 5.43 & 15.5 & 321.96 & 1679\\
    Profit&0.01 & 61.93 & 0.28 & 0.88 & 3.06 & 15.23 & 292.13 & 616\\
  	\hline
  \end{tabular}
  \caption{Summary statistics for Danish building fire losses data.}
  \label{tab:danishsummary}
\end{table}
\begin{table}[h!]
  \centering
  \begin{tabular}{lrrrrr}
  \hline
  Type & Distribution &   $LL_{is}$  &    $AIC_{is}$ &    $BIC_{is}$ &   $LL_{os}$ \\
  \hline
  Buildings  & IG & -2849.905 & 5703.811 & 5714.797 & -278.373 \\ 
        & GIG & -2848.965 & 5703.930 & 5720.410 & -277.869 \\ 
        & GG & -2675.870 & 5357.739 & 5374.219 & -263.656 \\
        & FIG & -2649.058 & 5306.117 & 5328.090 & -260.759 \\  
  \hline
  Content   & IG & -2100.065 & 4204.129 & 4214.757 & -213.877 \\ 
             & GIG & -2032.208 & 4070.415 & 4086.357 & -220.340 \\                 
                & GG & -1895.756 & 3797.512 & 3813.453 & -212.069 \\
                & FIG & -1879.272 & 3766.544 & 3787.799 & -209.191 \\
    \hline
    Profit      & IG & -327.101 & 658.203 & 666.815 & -18.264 \\ 
                & GIG & -312.282 & 630.564 & 643.483 & -18.597 \\ 
                & GG & -292.055 & 590.110 & 603.029 & -17.962 \\
                & FIG & -288.890 & 585.780 & 603.005 & -17.474 \\    
    \hline
    \end{tabular}
    \caption{In- and out-of-sample metrics of distributions fitted to Danish fire losses data.}
    \label{tab:danishmetrics}
\end{table}
The in-sample criterion and out-of-sample $LL_{os}$ are tabulated in Table \ref{tab:danishmetrics}. In this application, both the FIG is a clear favourite in all types of fire losses, since the in- and out-sample metrics coincide. That is, the in-sample metrics is lowest for the FIG as well as the highest for the out-of sample metric $LL_{os}$. It can therefore be concluded that for these loss data, a generalised model is more appropriate.
\section{Conclusions}\label{sec:con}
In this paper, we address the need for more flexible distributions without compromising on desirable distribution traits for positive data (Section \ref{sec:intro}). The method of Weibullisation is demonstrated (Section \ref{sec:genweibull}) and applied to the BTN distribution to yield the FIG distribution. The FIG has the desirable properties such as a low number of interpretable parameters, simple tractability, and finite moments. We provide many common statistical properties for using the FIG in practice. These are PDF, CDF, moments, MGF, and maximum likelihood estimation equations (Section \ref{sec:FIG} and \ref{sec:estim}). Regarding identifiability, a proof that the FIG parameters are identifiable, is provided (Section \ref{sec:identifia}). The applicability of the FIG is demonstrated on hand grip strength and insurance loss data where the FIG provides competitive fit in comparison to the IG and GIG distributions (Section \ref{sec:application}). Points for further research may include finite mixture modelling, kernel density smoothing, outlier detection, gamma regression, and reliability modelling.
\section*{Acknowledgments}
We would like to thank Prof. Angelo Spina from the Scientific High School E. Boggio Lera, Catania, Italy, for his assistance in determining the identifiability of the FIG distribution.
\subsection*{Financial disclosure}
 This work is based on the research supported in part by the National Research Foundation of South Africa Ref.: SRUG2204203965; RA171022270376, UID 119109; RA211204653274, Grant NO. 151035, as well as the Centre of Excellence in Mathematical and Statistical Sciences at the University of the Witwatersrand. The work of the author Muhammad Arashi, supported by the Iran National Science Foundation (INSF) Grant NO.4015320. Opinions expressed and conclusions arrived at are those of the authors and are not necessarily to be attributed to the funders of this work.
\subsection*{Conflict of interest}
The authors declare no potential conflict of interests.
%
\appendix
\section*{Lemmas}
\label{asec:A}
\subsubsection*{Lemma 1}
 Let $\alpha,\beta>0$. Then the following limit holds true
\begin{equation}
  \label{eq:lem1}
  \lim\limits_{x\to \infty} x^{k}\Gamma\left( \frac{\alpha}{\beta},x^\beta\right) =0 \text{\hspace{3mm}for }k\in\mathbb{R}.
 \end{equation}
If $k\le0$ both factors on the left-hand side of (\ref{eq:lem1}) tend to zero as $x$ tends to infinity.
 If $k>0$, by L'Hospital rule
 \begin{equation}
  \lim\limits_{x\to \infty} x^k\Gamma\left( \frac{\alpha}{\beta},x^\beta \right)
  =
  \lim\limits_{x\to \infty}
  \frac{x^{\alpha+\beta+k-1}}{\text{e}^{x^\beta}}\cdot
  \frac{\beta}{k}
  = 0\notag
 \end{equation}
 \subsubsection*{Lemma 2}
Let $\alpha,\beta>0$, then the following integral identity holds true
\begin{equation}
  \label{eq:lem2}
  \int_x^\infty t^r\Gamma\left( \frac{\alpha}{\beta},t^\beta\right)dt
  = \frac
  {\Gamma\left(\frac{\alpha+r+1}{\beta},x^\beta \right)
  -x^{r+1}\Gamma\left(\frac{\alpha}{\beta},x^\beta \right)}
  {r+1}.
 \end{equation}
 Let $y=t^\beta$, which implies $t=y^{\frac{1}{\beta}}$. Integrating by parts, where
 $v'(y)=\frac{r+1}{\beta}y^{\frac{r+1}{\beta}-1}$ and $u(y)=\Gamma\left({\alpha}/{\beta},y \right)$. The latter implies that $v(y)=y^{\frac{r+1}{\beta}}$ and $u'(y)=-y^{\frac{\alpha}{\beta}-1}\text{e}^{-y}$. The integral is evaluated as
\begin{gather}
  \int^{\infty}_{x}t^r\Gamma\left({\alpha}/{\beta},t^\beta \right)dt
  =
  (r+1)^{-1}
  \left. y^{\frac{r+1}{\beta}}\Gamma\left({\alpha}/{\beta},y \right)\right|^\infty_{x^\beta}
  -(r+1)^{-1}
  \int^{\infty}_{x^\beta}y^{\frac{r+1}{\beta}}\cdot
  \left( -y^{\frac{\alpha}{\beta}-1}\text{e}^{-y}\right) dy.\notag
 \end{gather}
 Noting from Lemma 1 that $\lim\limits_{x\to \infty} y^{\frac{r+1}{\beta}}\Gamma\left({\alpha}/{\beta},y \right)=0$
 the result follows.
\subsubsection*{Lemma 3}
Let $\alpha,\beta>0$. Then the following integral identity holds true
\begin{equation}
  \label{eq:lem3}
  \int_0^\infty x^r\Gamma\left( \frac{\alpha}{\beta},x^\beta\right)dx
  = \frac
  {\Gamma\left(\frac{\alpha+r+1}{\beta}\right)}
  {r+1}.
 \end{equation}
The result follows by evaluating the limit, $\lim_{t \to 0^+}$, over the integral in Lemma 2.

\section*{Declaration of Competing Interest}
The authors declare that they have no known competing
financial interests or personal relationships that could have
appeared to influence the work reported in this paper.

\nocite{*}

\bibliographystyle{elsarticle-harv} 
\bibliography{sample, GG}

\begin{thebibliography}{71}
\expandafter\ifx\csname natexlab\endcsname\relax\def\natexlab#1{#1}\fi
\providecommand{\url}[1]{\texttt{#1}}
\providecommand{\href}[2]{#2}
\providecommand{\path}[1]{#1}
\providecommand{\DOIprefix}{doi:}
\providecommand{\ArXivprefix}{arXiv:}
\providecommand{\URLprefix}{URL: }
\providecommand{\Pubmedprefix}{pmid:}
\providecommand{\doi}[1]{\href{http://dx.doi.org/#1}{\path{#1}}}
\providecommand{\Pubmed}[1]{\href{pmid:#1}{\path{#1}}}
\providecommand{\bibinfo}[2]{#2}
\ifx\xfnm\relax \def\xfnm[#1]{\unskip,\space#1}\fi
\bibitem[{Abid et~al.(2020)Abid, Al-Noor and Boshi}]{abid2020generalized}
\bibinfo{author}{Abid, S.H.}, \bibinfo{author}{Al-Noor, N.H.},
  \bibinfo{author}{Boshi, M.A.A.}, \bibinfo{year}{2020}.
\newblock \bibinfo{title}{{T}he {G}eneralized {G}amma-{E}xponentiated {W}eibull
  {D}istribution with its {P}roperties}.
\newblock \bibinfo{journal}{Al-Mustansiriyah Journal of Science}
  \bibinfo{volume}{31}, \bibinfo{pages}{30--37}.
\bibitem[{Agarwal and Kalla(1996)}]{IGGagarkall}
\bibinfo{author}{Agarwal, S.}, \bibinfo{author}{Kalla, S.},
  \bibinfo{year}{1996}.
\newblock \bibinfo{title}{{A} generalized gamma distribution and its
  application in reliabilty}.
\newblock \bibinfo{journal}{{C}ommunications in {S}tatistics-{T}heory and
  {M}ethods} \bibinfo{volume}{25}, \bibinfo{pages}{201--210}.
\newblock \DOIprefix\doi{10.1080/03610929608831688}.
\bibitem[{Agarwal and Al-Saleh(2001)}]{GGagar}
\bibinfo{author}{Agarwal, S.K.}, \bibinfo{author}{Al-Saleh, J.A.},
  \bibinfo{year}{2001}.
\newblock \bibinfo{title}{{G}eneralized gamma type distribution and its hazard
  rate function}.
\newblock \bibinfo{journal}{{C}ommunications in {S}tatistics-{T}heory and
  {M}ethods} \bibinfo{volume}{30}, \bibinfo{pages}{309--318}.
\newblock \DOIprefix\doi{10.1081/STA-100002033}.
\bibitem[{Akaike(1974)}]{akaike1974new}
\bibinfo{author}{Akaike, H.}, \bibinfo{year}{1974}.
\newblock \bibinfo{title}{A new look at the statistical model identification}.
\newblock \bibinfo{journal}{IEEE Transactions on Automatic Control}
  \bibinfo{volume}{19}, \bibinfo{pages}{716--723}.
\bibitem[{Ali et~al.(2008)Ali, Woo and Nadarajah}]{GGNadarajahapp}
\bibinfo{author}{Ali, M.}, \bibinfo{author}{Woo, J.},
  \bibinfo{author}{Nadarajah, S.}, \bibinfo{year}{2008}.
\newblock \bibinfo{title}{{G}eneralized gamma variables with drought
  application}.
\newblock \bibinfo{journal}{{J}ournal of the {K}orean {S}tatistical {S}ociety}
  \bibinfo{volume}{37}, \bibinfo{pages}{37--45}.
\newblock \DOIprefix\doi{10.1016/j.jkss.2007.09.002}.
\bibitem[{Almpanidis and Kotropoulos(2008)}]{GGTextapp}
\bibinfo{author}{Almpanidis, G.}, \bibinfo{author}{Kotropoulos, C.},
  \bibinfo{year}{2008}.
\newblock \bibinfo{title}{{P}honemic {S}egmentation using the {G}eneralised
  {G}amma {D}istribution and {S}mall {S}ample {B}ayesian {I}nformation
  {C}riterion}.
\newblock \bibinfo{journal}{{S}peech {C}ommunication} \bibinfo{volume}{50},
  \bibinfo{pages}{38--55}.
\newblock \DOIprefix\doi{10.1016/j.specom.2007.06.005}.
\bibitem[{Altun et~al.(2021)Altun, Korkmaz, El-Morshedy and Eliwa}]{2020NEG}
\bibinfo{author}{Altun, E.}, \bibinfo{author}{Korkmaz, M.C.},
  \bibinfo{author}{El-Morshedy, M.}, \bibinfo{author}{Eliwa, M.},
  \bibinfo{year}{2021}.
\newblock \bibinfo{title}{{T}he extended gamma distribution with regression
  model and applications}.
\newblock \bibinfo{journal}{AIMS Mathematics} \bibinfo{volume}{6}.
\bibitem[{Amoroso(1925)}]{Amoroso}
\bibinfo{author}{Amoroso, L.}, \bibinfo{year}{1925}.
\newblock \bibinfo{title}{{R}icerche intorno alla curva dei redditi}.
\newblock \bibinfo{journal}{{A}nnali di {M}atematica pura ed applicata}
  \bibinfo{volume}{2}, \bibinfo{pages}{123--159}.
\newblock \DOIprefix\doi{10.1007/BF02409935}.
\bibitem[{Bakery et~al.(2021)Bakery, Zakaria and Mohamed}]{2021DTGG}
\bibinfo{author}{Bakery, A.A.}, \bibinfo{author}{Zakaria, W.},
  \bibinfo{author}{Mohamed, O.M.K.S.K.}, \bibinfo{year}{2021}.
\newblock \bibinfo{title}{{A} {N}ew {D}ouble {T}runcated {G}eneralized {G}amma
  {M}odel with {S}ome {A}pplications}.
\newblock \bibinfo{journal}{Journal of Mathematics} \bibinfo{volume}{2021},
  \bibinfo{pages}{27}.
\newblock \URLprefix \url{https://doi.org/10.1155/2021/5500631},
  \DOIprefix\doi{10.1155/2021/5500631}.
\bibitem[{Barkauskas et~al.(2009)Barkauskas, Kronewitter, Lebrilla and
  Rocke}]{GGARMA}
\bibinfo{author}{Barkauskas, D.}, \bibinfo{author}{Kronewitter, S.},
  \bibinfo{author}{Lebrilla, C.}, \bibinfo{author}{Rocke, D.},
  \bibinfo{year}{2009}.
\newblock \bibinfo{title}{{A}nalysis of {MALDI FT-ICR} mass spectrometry data:
  {A} time series approach}.
\newblock \bibinfo{journal}{{A}nalytica {C}himica {A}cta}
  \bibinfo{volume}{648}, \bibinfo{pages}{207--214}.
\newblock \DOIprefix\doi{10.1016/j.aca.2009.06.064}.
\bibitem[{Barriga et~al.(2018)Barriga, Cordeiro, Dey, Cancho, Louzada and
  Suzuki}]{2018MOGG}
\bibinfo{author}{Barriga, G.}, \bibinfo{author}{Cordeiro, G.},
  \bibinfo{author}{Dey, D.}, \bibinfo{author}{Cancho, V.},
  \bibinfo{author}{Louzada, F.}, \bibinfo{author}{Suzuki, A.},
  \bibinfo{year}{2018}.
\newblock \bibinfo{title}{{T}he {M}arshall-{O}lkin generalized gamma
  distribution}.
\newblock \bibinfo{journal}{Communications for Statistical Applications and
  Methods} \bibinfo{volume}{25}, \bibinfo{pages}{245--261}.
\bibitem[{Bell(1988)}]{bellreparGG}
\bibinfo{author}{Bell, B.M.}, \bibinfo{year}{1988}.
\newblock \bibinfo{title}{Generalized gamma parameter estimation and moment
  evaluation}.
\newblock \bibinfo{journal}{Communications in Statistics-Theory and Methods}
  \bibinfo{volume}{17}, \bibinfo{pages}{507--517}.
\bibitem[{Bilankulu et~al.(2021)Bilankulu, Bekker and Marques}]{RGGBilankulu}
\bibinfo{author}{Bilankulu, V.}, \bibinfo{author}{Bekker, A.},
  \bibinfo{author}{Marques, F.}, \bibinfo{year}{2021}.
\newblock \bibinfo{title}{{T}he ratio of independent generalized gamma random
  variables with applications}.
\newblock \bibinfo{journal}{{C}omputational and {M}athematical {M}ethods}
  \bibinfo{volume}{3}, \bibinfo{pages}{e1061}.
\newblock \DOIprefix\doi{10.1002/cmm4.1061}.
\bibitem[{Bourguignon et~al.(2015)Bourguignon, do~Carmo Lima M.~S., Le{\~{a}}o,
  Nascimento, Pinho and Cordeiro}]{ENHTGBourguignon}
\bibinfo{author}{Bourguignon, M.}, \bibinfo{author}{do~Carmo Lima M.~S.},
  \bibinfo{author}{Le{\~{a}}o, J.}, \bibinfo{author}{Nascimento, A.D.C.},
  \bibinfo{author}{Pinho, L.}, \bibinfo{author}{Cordeiro, G.},
  \bibinfo{year}{2015}.
\newblock \bibinfo{title}{{A} new generalized gamma distribution with
  applications}.
\newblock \bibinfo{journal}{{A}merican {J}ournal of {M}athematical and
  {M}anagement {S}ciences} \bibinfo{volume}{34}, \bibinfo{pages}{309--342}.
\newblock \DOIprefix\doi{10.1080/01966324.2015.1040178}.
\bibitem[{Chen et~al.(2011)Chen, Karagiannidis, Lu and Cao}]{GGGwirelessapp}
\bibinfo{author}{Chen, Y.}, \bibinfo{author}{Karagiannidis, G.},
  \bibinfo{author}{Lu, H.}, \bibinfo{author}{Cao, N.}, \bibinfo{year}{2011}.
\newblock \bibinfo{title}{{N}ovel {A}pproximations to the {S}tatistics of
  {P}roducts of {I}ndependent {R}andom {V}ariables and {T}heir {A}pplications
  in {W}ireless {C}ommunications}.
\newblock \bibinfo{journal}{{IEEE} {T}ransactions on {V}ehicular {T}echnology}
  \bibinfo{volume}{61}, \bibinfo{pages}{443--454}.
\newblock \DOIprefix\doi{10.1109/TVT.2011.2178441}.
\bibitem[{Consul and Jain(1971)}]{GGCosulJain}
\bibinfo{author}{Consul, P.C.}, \bibinfo{author}{Jain, G.C.},
  \bibinfo{year}{1971}.
\newblock \bibinfo{title}{{O}n the log-gamma distribution and its properties}.
\newblock \bibinfo{journal}{{S}tatistische {H}efte} \bibinfo{volume}{12},
  \bibinfo{pages}{100--106}.
\newblock \DOIprefix\doi{10.1007/BF02922944}.
\bibitem[{Cordeiro et~al.(2012)Cordeiro, Castellares, Montenegro and
  de~Castro}]{BGG}
\bibinfo{author}{Cordeiro, G.M.}, \bibinfo{author}{Castellares, F.},
  \bibinfo{author}{Montenegro, L.C.}, \bibinfo{author}{de~Castro, M.},
  \bibinfo{year}{2012}.
\newblock \bibinfo{title}{The beta generalized gamma distribution}.
\newblock \bibinfo{journal}{Statistics} \bibinfo{volume}{47},
  \bibinfo{pages}{888--900}.
\bibitem[{Cordeiro et~al.(2011)Cordeiro, Ortega and Silva}]{EGG}
\bibinfo{author}{Cordeiro, G.M.}, \bibinfo{author}{Ortega, E.M.},
  \bibinfo{author}{Silva, G.O.}, \bibinfo{year}{2011}.
\newblock \bibinfo{title}{{T}he exponentiated generalized gamma distribution
  with application to lifetime data}.
\newblock \bibinfo{journal}{{J}ournal of {S}tatistical {C}omputation and
  {S}imulation} \bibinfo{volume}{81}, \bibinfo{pages}{827--842}.
\newblock \DOIprefix\doi{10.1080/00949650903517874}.
\bibitem[{Cordeiro et~al.(2014)Cordeiro, Pescim and Ortega}]{KBGG}
\bibinfo{author}{Cordeiro, G.M.}, \bibinfo{author}{Pescim, R.},
  \bibinfo{author}{Ortega, E.}, \bibinfo{year}{2014}.
\newblock \bibinfo{title}{The {K}ummer {B}eta {G}eneralized {G}amma
  {D}istribution}.
\newblock \bibinfo{journal}{{J}ournal of {D}ata {S}cience}
  \bibinfo{volume}{12}, \bibinfo{pages}{661--697}.
\newblock \DOIprefix\doi{10.6339/JDS.201410_12(4).0006}.
\bibitem[{Cover and Tomas(1999)}]{cover1999elements}
\bibinfo{author}{Cover, T.M.}, \bibinfo{author}{Tomas, J.A.},
  \bibinfo{year}{1999}.
\newblock \bibinfo{title}{{E}lements of information theory}.
\newblock \bibinfo{publisher}{{J}ohn {W}iley \& {S}ons}.
\bibitem[{De~Pascoa et~al.(2011)De~Pascoa, Ortega and Cordeiro}]{KGG}
\bibinfo{author}{De~Pascoa, M.A.R.}, \bibinfo{author}{Ortega, E.M.},
  \bibinfo{author}{Cordeiro, G.M.}, \bibinfo{year}{2011}.
\newblock \bibinfo{title}{The kumaraswamy generalized gamma distribution with
  application in survival analysis}.
\newblock \bibinfo{journal}{Statistical methodology} \bibinfo{volume}{8},
  \bibinfo{pages}{411--433}.
\bibitem[{Du and Chen(2016)}]{GGasymtotics}
\bibinfo{author}{Du, L.}, \bibinfo{author}{Chen, S.}, \bibinfo{year}{2016}.
\newblock \bibinfo{title}{{A}symptotic properties for distributions and
  densities of extremes from generalized gamma distribution}.
\newblock \bibinfo{journal}{{J}ournal of the {K}orean {S}tatistical {S}ociety}
  \bibinfo{volume}{45}, \bibinfo{pages}{188--198}.
\newblock \DOIprefix\doi{10.1016/j.jkss.2015.09.005}.
\bibitem[{Folks and Chhikara(1978)}]{folks1978inverse}
\bibinfo{author}{Folks, J.}, \bibinfo{author}{Chhikara, R.},
  \bibinfo{year}{1978}.
\newblock \bibinfo{title}{The inverse gaussian distribution and its statistical
  application—a review}.
\newblock \bibinfo{journal}{Journal of the Royal Statistical Society: Series B
  (Methodological)} \bibinfo{volume}{40}, \bibinfo{pages}{263--275}.
\bibitem[{Gradshteyn and Ryhzhik(2007)}]{Gr<3dstheyn}
\bibinfo{author}{Gradshteyn, I.}, \bibinfo{author}{Ryhzhik, I.},
  \bibinfo{year}{2007}.
\newblock \bibinfo{title}{Table of Integrals, Series and Products}.
\newblock \bibinfo{publisher}{Academic Press:Cambridge, MA, USA}.
\bibitem[{Harris et~al.(2020)Harris, Millman, van~der Walt, Gommers, Virtanen,
  Cournapeau, Wieser, Taylor, Berg, Smith, Kern, Picus, Hoyer, van Kerkwijk,
  Brett, Haldane, Fernández~del Río, Wiebe, Peterson, Gérard-Marchant,
  Sheppard, Reddy, Weckesser, Abbasi, Gohlke and Oliphant}]{numpy}
\bibinfo{author}{Harris, C.R.}, \bibinfo{author}{Millman, K.J.},
  \bibinfo{author}{van~der Walt, S.J.}, \bibinfo{author}{Gommers, R.},
  \bibinfo{author}{Virtanen, P.}, \bibinfo{author}{Cournapeau, D.},
  \bibinfo{author}{Wieser, E.}, \bibinfo{author}{Taylor, J.},
  \bibinfo{author}{Berg, S.}, \bibinfo{author}{Smith, N.J.},
  \bibinfo{author}{Kern, R.}, \bibinfo{author}{Picus, M.},
  \bibinfo{author}{Hoyer, S.}, \bibinfo{author}{van Kerkwijk, M.H.},
  \bibinfo{author}{Brett, M.}, \bibinfo{author}{Haldane, A.},
  \bibinfo{author}{Fernández~del Río, J.}, \bibinfo{author}{Wiebe, M.},
  \bibinfo{author}{Peterson, P.}, \bibinfo{author}{Gérard-Marchant, P.},
  \bibinfo{author}{Sheppard, K.}, \bibinfo{author}{Reddy, T.},
  \bibinfo{author}{Weckesser, W.}, \bibinfo{author}{Abbasi, H.},
  \bibinfo{author}{Gohlke, C.}, \bibinfo{author}{Oliphant, T.E.},
  \bibinfo{year}{2020}.
\newblock \bibinfo{title}{Array programming with {NumPy}}.
\newblock \bibinfo{journal}{Nature} \bibinfo{volume}{585},
  \bibinfo{pages}{357–362}.
\newblock \DOIprefix\doi{10.1038/s41586-020-2649-2}.
\bibitem[{Johansson(2018)}]{mpmath}
\bibinfo{author}{Johansson, F.}, \bibinfo{year}{2018}.
\newblock \bibinfo{title}{mpmath: a {P}ython library for arbitrary-precision
  floating-point arithmetic (version 1.1.0)}.
\bibitem[{Johnson and Kotz(1972)}]{johnson1972power}
\bibinfo{author}{Johnson, N.}, \bibinfo{author}{Kotz, S.},
  \bibinfo{year}{1972}.
\newblock \bibinfo{title}{Power transformations of gamma variables}.
\newblock \bibinfo{journal}{Biometrika} , \bibinfo{pages}{226--229}.
\bibitem[{Johnson~N.L. et~al.(1970)Johnson~N.L., Kotz and
  Johnson}]{johnson1970continuous}
\bibinfo{author}{Johnson~N.L., N.}, \bibinfo{author}{Kotz, S.},
  \bibinfo{author}{Johnson, N.L.}, \bibinfo{year}{1970}.
\newblock \bibinfo{title}{Continuous univariate distributions}.
\newblock \bibinfo{journal}{AGRIS} .
\bibitem[{Jones(2015)}]{jones15}
\bibinfo{author}{Jones, M.}, \bibinfo{year}{2015}.
\newblock \bibinfo{title}{On {F}amilies of {D}istributions with {S}hape
  {P}arameters}.
\newblock \bibinfo{journal}{International {S}tatistical {R}eview}
  \bibinfo{volume}{83}, \bibinfo{pages}{175--192}.
\bibitem[{Jorgensen(2012)}]{GIGbook}
\bibinfo{author}{Jorgensen, B.}, \bibinfo{year}{2012}.
\newblock \bibinfo{title}{Statistical properties of the generalized inverse
  Gaussian distribution}. volume~\bibinfo{volume}{9}.
\newblock \bibinfo{publisher}{Springer Science \& Business Media}.
\bibitem[{Kalla and Al-Saqabi(2001)}]{GIGGKall}
\bibinfo{author}{Kalla, S.}, \bibinfo{author}{Al-Saqabi, B.},
  \bibinfo{year}{2001}.
\newblock \bibinfo{title}{{F}urther results on a unified form of gamma-type
  distributions}.
\newblock \bibinfo{journal}{{F}ractional {C}alculus and {A}pplied {A}nalysis}
  \bibinfo{volume}{4}, \bibinfo{pages}{91--100}.
\bibitem[{Kaneko(2003)}]{CoaleNuptiDemo}
\bibinfo{author}{Kaneko, R.}, \bibinfo{year}{2003}.
\newblock \bibinfo{title}{Elaboration of the coale-mcneil nuptiality model as
  the generalized log gamma distribution: A new identity and empirical
  enhancements}.
\newblock \bibinfo{journal}{Demographic Research} \bibinfo{volume}{9},
  \bibinfo{pages}{223--262}.
\bibitem[{Kaniadakis(2021)}]{2021kappaGG}
\bibinfo{author}{Kaniadakis, G.}, \bibinfo{year}{2021}.
\newblock \bibinfo{title}{New power-law tailed distributions emerging in
  $\kappa$-statistics}.
\newblock \bibinfo{journal}{Europhysics Letters} \bibinfo{volume}{133},
  \bibinfo{pages}{10002}.
\bibitem[{Kleiber and Kotz(2003)}]{kleiber2003statistical}
\bibinfo{author}{Kleiber, C.}, \bibinfo{author}{Kotz, S.},
  \bibinfo{year}{2003}.
\newblock \bibinfo{title}{{S}tatistical {S}ize {D}istributions in {E}conomics
  and {A}ctuarial {S}ciences}. volume \bibinfo{volume}{470}.
\newblock \bibinfo{publisher}{{J}ohn {W}iley \& {S}ons}.
\bibitem[{Kobayashi(1991a)}]{diffractionGG}
\bibinfo{author}{Kobayashi, K.}, \bibinfo{year}{1991}a.
\newblock \bibinfo{title}{{O}n {G}eneralized {G}amma {F}unctions {O}ccurring in
  {D}iffraction {T}heory}.
\newblock \bibinfo{journal}{{J}ournal of the {P}hysical {S}ociety of {J}apan}
  \bibinfo{volume}{60}, \bibinfo{pages}{1501--1512}.
\newblock \DOIprefix\doi{10.1143/JPSJ.60.1501}.
\bibitem[{Kobayashi(1991b)}]{kobayashiGG}
\bibinfo{author}{Kobayashi, K.}, \bibinfo{year}{1991}b.
\newblock \bibinfo{title}{On generalized gamma functions occurring in
  diffraction theory}.
\newblock \bibinfo{journal}{Journal of the Physical Society of Japan}
  \bibinfo{volume}{60}, \bibinfo{pages}{1501--1512}.
\bibitem[{Lee and Gross(1991)}]{GGleegros}
\bibinfo{author}{Lee, M.T.}, \bibinfo{author}{Gross, A.J.},
  \bibinfo{year}{1991}.
\newblock \bibinfo{title}{{L}ifetime distributions under unknown environment}.
\newblock \bibinfo{journal}{{J}ournal of {S}tatistical {P}lanning and
  {I}nference} \bibinfo{volume}{29}, \bibinfo{pages}{137--143}.
\newblock \DOIprefix\doi{10.1016/0378-3758(92)90128-F}.
\bibitem[{Lees-Miller et~al.(2010)Lees-Miller, Hammersley and
  Wilson}]{lees2010theoretical}
\bibinfo{author}{Lees-Miller, J.}, \bibinfo{author}{Hammersley, J.},
  \bibinfo{author}{Wilson, R.}, \bibinfo{year}{2010}.
\newblock \bibinfo{title}{Theoretical maximum capacity as benchmark for empty
  vehicle redistribution in personal rapid transit}.
\newblock \bibinfo{journal}{Transportation Research Record: Journal of the
  Transportation Research Board} , \bibinfo{pages}{76--83}.
\bibitem[{Lehmann and Casella(2006)}]{GGLehman}
\bibinfo{author}{Lehmann, E.L.}, \bibinfo{author}{Casella, G.},
  \bibinfo{year}{2006}.
\newblock \bibinfo{title}{{T}heory of {P}oint {E}stimation}.
\newblock \bibinfo{publisher}{{S}pringer {S}cience \& {B}usiness {M}edia}.
\bibitem[{Ley(2015)}]{ley14flex}
\bibinfo{author}{Ley, C.}, \bibinfo{year}{2015}.
\newblock \bibinfo{title}{Flexible modelling in statistics: past, present and
  future}.
\newblock \bibinfo{journal}{Journal de la Societe Francaise de Statistique}
  \bibinfo{volume}{156}.
\bibitem[{Ley et~al.(2021)Ley, Babi{\'c} and Craens}]{ley2021flexible}
\bibinfo{author}{Ley, C.}, \bibinfo{author}{Babi{\'c}, S.},
  \bibinfo{author}{Craens, D.}, \bibinfo{year}{2021}.
\newblock \bibinfo{title}{Flexible models for complex data with applications}.
\newblock \bibinfo{journal}{Annual Review of Statistics and Its Application}
  \bibinfo{volume}{8}, \bibinfo{pages}{369--391}.
\bibitem[{Lucena et~al.(2015)Lucena, Herm and Cordeiro}]{transmGG}
\bibinfo{author}{Lucena, S.E.}, \bibinfo{author}{Herm, A.},
  \bibinfo{author}{Cordeiro, G.M.}, \bibinfo{year}{2015}.
\newblock \bibinfo{title}{The transmuted generalized gamma distribution:
  Properties and application}.
\newblock \bibinfo{journal}{Journal of Data Science} \bibinfo{volume}{13},
  \bibinfo{pages}{187--206}.
\bibitem[{Malhotra et~al.(2009)Malhotra, Sharma and Kaler}]{GGfading}
\bibinfo{author}{Malhotra, J.}, \bibinfo{author}{Sharma, A.K.},
  \bibinfo{author}{Kaler, R.}, \bibinfo{year}{2009}.
\newblock \bibinfo{title}{{O}n the performance analysis of wireless receiver
  using generalized-gamma fading model}.
\newblock \bibinfo{journal}{annals of telecommunications-annales des
  t{\'e}l{\'e}communications} \bibinfo{volume}{64}, \bibinfo{pages}{147--153}.
\bibitem[{Malik(1967)}]{QuoGG}
\bibinfo{author}{Malik, H.J.}, \bibinfo{year}{1967}.
\newblock \bibinfo{title}{{E}xact {D}istribution of the {Q}uotient of
  {I}ndependent {G}eneralized {G}amma {V}ariables}.
\newblock \bibinfo{journal}{{C}anadian {M}athematical {B}ulletin}
  \bibinfo{volume}{10}, \bibinfo{pages}{463--465}.
\newblock \DOIprefix\doi{10.4153/CMB-1967-045-7}.
\bibitem[{McDonald and Xu(1995)}]{Gbeta}
\bibinfo{author}{McDonald, J.}, \bibinfo{author}{Xu, Y.}, \bibinfo{year}{1995}.
\newblock \bibinfo{title}{{A} generalization of the beta distribution with
  applications}.
\newblock \bibinfo{journal}{{J}ournal of {E}conometrics} \bibinfo{volume}{66},
  \bibinfo{pages}{133--152}.
\newblock \DOIprefix\doi{10.1016/0304-4076(94)01612-4}.
\bibitem[{McLeish(1982)}]{mcleish1982robust}
\bibinfo{author}{McLeish, D.L.}, \bibinfo{year}{1982}.
\newblock \bibinfo{title}{A {R}obust {A}lternative to the {N}normal
  {D}istribution}.
\newblock \bibinfo{journal}{The {C}anadian {J}ournal of {S}tatistics} ,
  \bibinfo{pages}{89--102}.
\bibitem[{Mead(2015)}]{GIGGMead}
\bibinfo{author}{Mead, M.}, \bibinfo{year}{2015}.
\newblock \bibinfo{title}{{G}eneralized inverse gamma distribution and its
  application in reliability}.
\newblock \bibinfo{journal}{{C}ommunications in {S}tatistics-{T}heory and
  {M}ethods} \bibinfo{volume}{44}, \bibinfo{pages}{1426--1435}.
\newblock \DOIprefix\doi{10.1080/03610926.2013.768667}.
\bibitem[{Mead et~al.(2018)Mead, Nassar and Dey}]{2018MGG}
\bibinfo{author}{Mead, M.}, \bibinfo{author}{Nassar, M.}, \bibinfo{author}{Dey,
  S.}, \bibinfo{year}{2018}.
\newblock \bibinfo{title}{{A} generalization of generalized gamma
  distributions}.
\newblock \bibinfo{journal}{Pakistan Journal of Statistics and Operation
  Research} \bibinfo{volume}{14}, \bibinfo{pages}{121--138}.
\bibitem[{Nadarajah and Haghighi(2011)}]{SaraGE}
\bibinfo{author}{Nadarajah, S.}, \bibinfo{author}{Haghighi, F.},
  \bibinfo{year}{2011}.
\newblock \bibinfo{title}{{A}n extension of the exponential distribution}.
\newblock \bibinfo{journal}{{S}tatistics} \bibinfo{volume}{45},
  \bibinfo{pages}{543--558}.
\newblock \DOIprefix\doi{10.1080/02331881003678678}.
\bibitem[{Nadarajahah and G.(2007)}]{GGNaGup}
\bibinfo{author}{Nadarajahah, S.}, \bibinfo{author}{G., A.K.},
  \bibinfo{year}{2007}.
\newblock \bibinfo{title}{A generalized gamma distribution with application to
  drought data}.
\newblock \bibinfo{journal}{{M}athematics and {C}omputers in {S}imulation}
  \bibinfo{volume}{74}, \bibinfo{pages}{1--7}.
\newblock \DOIprefix\doi{10.1016/j.matcom.2006.04.004}.
\bibitem[{Pauw et~al.(2010)Pauw, Bekker and Roux}]{pauw2010densities}
\bibinfo{author}{Pauw, J.}, \bibinfo{author}{Bekker, A.},
  \bibinfo{author}{Roux, J.}, \bibinfo{year}{2010}.
\newblock \bibinfo{title}{{D}ensities of composite {W}eibullized generalized
  gamma variables: theory and methods}.
\newblock \bibinfo{journal}{{S}outh {A}frican {S}tatistical {J}ournal}
  \bibinfo{volume}{44}, \bibinfo{pages}{17--42}.
\bibitem[{Pham-Gia and Duong(1989)}]{PHAMGIA1989}
\bibinfo{author}{Pham-Gia, T.}, \bibinfo{author}{Duong, Q.},
  \bibinfo{year}{1989}.
\newblock \bibinfo{title}{The generalized beta- and f-distributions in
  statistical modelling}.
\newblock \bibinfo{journal}{Mathematical and Computer Modelling}
  \bibinfo{volume}{12}, \bibinfo{pages}{1613--1625}.
\newblock \DOIprefix\doi{https://doi.org/10.1016/0895-7177(89)90337-3}.
\bibitem[{Prentice(1974)}]{prenticelogGG}
\bibinfo{author}{Prentice, R.L.}, \bibinfo{year}{1974}.
\newblock \bibinfo{title}{A log gamma model and its maximum likelihood
  estimation}.
\newblock \bibinfo{journal}{Biometrika} \bibinfo{volume}{61},
  \bibinfo{pages}{539--544}.
\bibitem[{Priyadarshani and Oluyede(2015)}]{priyadarshaniWGG}
\bibinfo{author}{Priyadarshani, H.A.}, \bibinfo{author}{Oluyede, B.O.},
  \bibinfo{year}{2015}.
\newblock \bibinfo{title}{Theoretical properties of the weighted generalized
  gamma and related distributions}.
\newblock \bibinfo{journal}{Probability in the Engineering and Informational
  Sciences} \bibinfo{volume}{29}, \bibinfo{pages}{421--432}.
\bibitem[{Punzo and Bagnato(2021)}]{punzo2021multivariate}
\bibinfo{author}{Punzo, A.}, \bibinfo{author}{Bagnato, L.},
  \bibinfo{year}{2021}.
\newblock \bibinfo{title}{The multivariate tail-inflated normal distribution
  and its application in finance}.
\newblock \bibinfo{journal}{Journal of {S}tatistical {C}omputation and
  {S}imulation} \bibinfo{volume}{91}, \bibinfo{pages}{1--36}.
\bibitem[{Ramos et~al.(2021)Ramos, Mota, Ferreira, Ramos, Tomazella and
  Louzanda}]{IGGBayesian}
\bibinfo{author}{Ramos, P.}, \bibinfo{author}{Mota, A.},
  \bibinfo{author}{Ferreira, P.}, \bibinfo{author}{Ramos, E.},
  \bibinfo{author}{Tomazella, V.}, \bibinfo{author}{Louzanda, F.},
  \bibinfo{year}{2021}.
\newblock \bibinfo{title}{{B}ayesian analysis of the inverse generalized gamma
  distribution using objective priors}.
\newblock \bibinfo{journal}{{J}ournal of {S}tatistical {C}omputation and
  {S}imulation} \bibinfo{volume}{91}, \bibinfo{pages}{786--816}.
\newblock \DOIprefix\doi{10.1080/00949655.2020.1830991}.
\bibitem[{Rigby and Stasinopoulos(2005)}]{gamlss}
\bibinfo{author}{Rigby, R.A.}, \bibinfo{author}{Stasinopoulos, D.M.},
  \bibinfo{year}{2005}.
\newblock \bibinfo{title}{Generalized additive models for location, scale and
  shape,(with discussion)}.
\newblock \bibinfo{journal}{Applied Statistics} \bibinfo{volume}{54},
  \bibinfo{pages}{507--554}.
\bibitem[{Saieed et~al.(2020)Saieed, Abdulla and Hayawi}]{2020IGG}
\bibinfo{author}{Saieed, H.}, \bibinfo{author}{Abdulla, M.},
  \bibinfo{author}{Hayawi, H.}, \bibinfo{year}{2020}.
\newblock \bibinfo{title}{Inverse generalized gamma distribution with it's
  properties}.
\newblock \bibinfo{journal}{IRAQI JOURNAL OF STATISTICAL SCIENCES}
  \bibinfo{volume}{17}, \bibinfo{pages}{64--71}.
\bibitem[{Samuelson et~al.(2006)Samuelson, D. and E.}]{IGbook}
\bibinfo{author}{Samuelson, P.A.}, \bibinfo{author}{D., M.},
  \bibinfo{author}{E., A.}, \bibinfo{year}{2006}.
\newblock \bibinfo{title}{Louis Bachelier's theory of speculation: the origins
  of modern finance}.
\newblock \bibinfo{publisher}{JSTOR}.
\bibitem[{Schwarz(1978)}]{schwarz1978estimating}
\bibinfo{author}{Schwarz, G.}, \bibinfo{year}{1978}.
\newblock \bibinfo{title}{Estimating the dimension of a model}.
\newblock \bibinfo{journal}{The Annals of Statistics} \bibinfo{volume}{6},
  \bibinfo{pages}{461--464}.
\bibitem[{Seshadri(2012)}]{seshadri2012inverse}
\bibinfo{author}{Seshadri, V.}, \bibinfo{year}{2012}.
\newblock \bibinfo{title}{The Inverse Gaussian Distribution: Statistical Theory
  and Applications}. volume \bibinfo{volume}{137}.
\newblock \bibinfo{publisher}{Springer Science \& Business Media}.
\bibitem[{Stacy and Hoshkin(1962)}]{ggStacy}
\bibinfo{author}{Stacy, E.W.}, \bibinfo{author}{Hoshkin, N.},
  \bibinfo{year}{1962}.
\newblock \bibinfo{title}{{A} {G}eneralization of the {G}amma {D}istribution}.
\newblock \bibinfo{journal}{The {A}nnals of {M}athematical {S}tatistics}
  \bibinfo{volume}{33}, \bibinfo{pages}{1187--1192}.
\bibitem[{Subbotin(1923)}]{subbotin23}
\bibinfo{author}{Subbotin, M.T.}, \bibinfo{year}{1923}.
\newblock \bibinfo{title}{On the {L}aw of {F}requency of {E}rror}.
\newblock \bibinfo{journal}{Mathematicheskii {S}bornik} \bibinfo{volume}{31},
  \bibinfo{pages}{296--301}.
\bibitem[{Tweedie(1957)}]{tweedie1957statistical}
\bibinfo{author}{Tweedie, M.}, \bibinfo{year}{1957}.
\newblock \bibinfo{title}{Statistical properties of inverse gaussian
  distributions. i}.
\newblock \bibinfo{journal}{The Annals of Mathematical Statistics}
  \bibinfo{volume}{28}, \bibinfo{pages}{362--377}.
\bibitem[{Umar and Yahya(2021)}]{2021NEGG}
\bibinfo{author}{Umar, M.}, \bibinfo{author}{Yahya, W.}, \bibinfo{year}{2021}.
\newblock \bibinfo{title}{{A} {N}ew {E}xponential-{G}amma {D}istribution with
  {A}pplications}.
\newblock \bibinfo{journal}{Journal of Modern Applied Statistical Methods} .
\bibitem[{Vallejos et~al.(2019)Vallejos, Ormaz{\'a}bal, Borotto and
  Astudillo}]{2018kappaG}
\bibinfo{author}{Vallejos, A.}, \bibinfo{author}{Ormaz{\'a}bal, I.},
  \bibinfo{author}{Borotto, F.A.}, \bibinfo{author}{Astudillo, H.F.},
  \bibinfo{year}{2019}.
\newblock \bibinfo{title}{{A} new $\kappa$-deformed parametric model for the
  size distribution of wealth}.
\newblock \bibinfo{journal}{Physica A: Statistical Mechanics and its
  Applications} \bibinfo{volume}{514}, \bibinfo{pages}{819--829}.
\bibitem[{Van~Niekerk et~al.(2017)Van~Niekerk, Bekker and
  Arashi}]{van2017gamma}
\bibinfo{author}{Van~Niekerk, J.}, \bibinfo{author}{Bekker, A.},
  \bibinfo{author}{Arashi, M.}, \bibinfo{year}{2017}.
\newblock \bibinfo{title}{{A} gamma-mixture class of distributions with
  {B}ayesian application}.
\newblock \bibinfo{journal}{{C}ommunications in {S}tatistics-{S}imulation and
  {C}omputation} \bibinfo{volume}{46}, \bibinfo{pages}{8152--8165}.
\newblock \DOIprefix\doi{10.1080/03610918.2016.1267754}.
\bibitem[{Virtanen et~al.(2020)Virtanen, Gommers, Oliphant, Haberland, Reddy,
  Cournapeau, Burovski, Peterson, Weckesser, Bright, {van der Walt}, Brett,
  Wilson, Millman, Mayorov, Nelson, Jones, Kern, Larson, Carey, Polat, Feng,
  Moore, {VanderPlas}, Laxalde, Perktold, Cimrman, Henriksen, Quintero, Harris,
  Archibald, Ribeiro, Pedregosa, {van Mulbregt} and {SciPy 1.0
  Contributors}}]{scipy}
\bibinfo{author}{Virtanen, P.}, \bibinfo{author}{Gommers, R.},
  \bibinfo{author}{Oliphant, T.E.}, \bibinfo{author}{Haberland, M.},
  \bibinfo{author}{Reddy, T.}, \bibinfo{author}{Cournapeau, D.},
  \bibinfo{author}{Burovski, E.}, \bibinfo{author}{Peterson, P.},
  \bibinfo{author}{Weckesser, W.}, \bibinfo{author}{Bright, J.},
  \bibinfo{author}{{van der Walt}, S.J.}, \bibinfo{author}{Brett, M.},
  \bibinfo{author}{Wilson, J.}, \bibinfo{author}{Millman, K.J.},
  \bibinfo{author}{Mayorov, N.}, \bibinfo{author}{Nelson, A.R.J.},
  \bibinfo{author}{Jones, E.}, \bibinfo{author}{Kern, R.},
  \bibinfo{author}{Larson, E.}, \bibinfo{author}{Carey, C.J.},
  \bibinfo{author}{Polat, {\.I}.}, \bibinfo{author}{Feng, Y.},
  \bibinfo{author}{Moore, E.W.}, \bibinfo{author}{{VanderPlas}, J.},
  \bibinfo{author}{Laxalde, D.}, \bibinfo{author}{Perktold, J.},
  \bibinfo{author}{Cimrman, R.}, \bibinfo{author}{Henriksen, I.},
  \bibinfo{author}{Quintero, E.A.}, \bibinfo{author}{Harris, C.R.},
  \bibinfo{author}{Archibald, A.M.}, \bibinfo{author}{Ribeiro, A.H.},
  \bibinfo{author}{Pedregosa, F.}, \bibinfo{author}{{van Mulbregt}, P.},
  \bibinfo{author}{{SciPy 1.0 Contributors}}, \bibinfo{year}{2020}.
\newblock \bibinfo{title}{{{SciPy} 1.0: Fundamental Algorithms for Scientific
  Computing in Python}}.
\newblock \bibinfo{journal}{Nature Methods} \bibinfo{volume}{17},
  \bibinfo{pages}{261--272}.
\newblock \DOIprefix\doi{10.1038/s41592-019-0686-2}.
\bibitem[{Wagener et~al.(2021)Wagener, Bekker and Arashi}]{BTGN}
\bibinfo{author}{Wagener, M.}, \bibinfo{author}{Bekker, A.},
  \bibinfo{author}{Arashi, M.}, \bibinfo{year}{2021}.
\newblock \bibinfo{title}{Mastering the body and tail shape of a distribution}.
\newblock \bibinfo{journal}{Mathematics} \bibinfo{volume}{9},
  \bibinfo{pages}{2648}.
\bibitem[{Yadav(2019)}]{2019IEG}
\bibinfo{author}{Yadav, A.S.}, \bibinfo{year}{2019}.
\newblock \bibinfo{title}{The inverted exponentiated gamma distribution: A
  heavy-tailed model with upside down bathtubshaped hazard rate}.
\newblock \bibinfo{journal}{Statistica} \bibinfo{volume}{79},
  \bibinfo{pages}{339--360}.
\bibitem[{Zografos and Balakrishnan(2009)}]{NHGBala}
\bibinfo{author}{Zografos, K.}, \bibinfo{author}{Balakrishnan, N.},
  \bibinfo{year}{2009}.
\newblock \bibinfo{title}{{O}n families of beta and generalized gamma-generated
  distributions and associated inference}.
\newblock \bibinfo{journal}{{S}tatistical {M}ethodology} \bibinfo{volume}{6},
  \bibinfo{pages}{344--362}.
\newblock \DOIprefix\doi{10.1016/j.stamet.2008.12.003}.

\end{thebibliography}

\end{document}